\documentclass[10pt,a4paper,twoside]{amsart}
\usepackage[bookmarks=false]{hyperref}
\hypersetup{colorlinks=true,allcolors = blue}
\usepackage{amsmath}
\usepackage{amsfonts}
\usepackage{amssymb}
\usepackage[applemac]{inputenc}
\usepackage[english]{babel}
\usepackage{verbatim}
\usepackage{graphicx}
\usepackage[all]{xy}
\xyoption{poly}
\usepackage{enumitem}
\usepackage{multirow}
\usepackage{longtable}
\usepackage{booktabs}
\usepackage{stmaryrd}
\usepackage[usenames,dvipsnames]{color}
\numberwithin{equation}{section}
	\usepackage{amsthm}
	\usepackage{aliascnt}
	\theoremstyle{plain}
	
	\newaliascnt{thm}{equation}
	\newtheorem{main}{Theorem}
	
	\newaliascnt{main_cor}{main}
	
	\aliascntresetthe{main_cor}
	\newaliascnt{prop}{thm}
	\newtheorem{prop}[prop]{Proposition}
	\aliascntresetthe{prop}
	\newaliascnt{lem}{thm}
	\newtheorem{lem}[lem]{Lemma}
	\aliascntresetthe{lem}

	\newaliascnt{cor}{thm}
	\newtheorem{cor}[cor]{Corollary}
	\aliascntresetthe{cor}
	
	\newaliascnt{conj}{thm}
	
	\aliascntresetthe{conj}
	\theoremstyle{definition}
	\newaliascnt{defn}{thm}
	\newtheorem{defn}[defn]{Definition}
	\aliascntresetthe{defn}
	\theoremstyle{remark}
	\newaliascnt{notation}{thm}
	\newtheorem{notation}[notation]{Notation}
	\aliascntresetthe{notation}
	\newaliascnt{rem}{thm}
	\newtheorem{rem}[rem]{Remark}
	\aliascntresetthe{rem}
	\newaliascnt{ass}{thm}
	
	\aliascntresetthe{ass}
	\newaliascnt{exmp}{thm}
	\newtheorem{exmp}[exmp]{Example}
	\aliascntresetthe{exmp}
\usepackage{array}
\usepackage[capitalise]{cleveref}
	\Crefname{equation}{}{}
	\crefname{equation}{}{}
	\Crefname{thm}{Theorem}{Theorems}
	\Crefname{main}{Theorem}{Theorems}
	\Crefname{main_cor}{Corollary}{Corollaries}
	\Crefname{prop}{Proposition}{Propositions}
	\Crefname{lem}{Lemma}{Lemmata}
	\Crefname{cor}{Corollary}{Corollaries}
	\Crefname{defn}{Definition}{Definitions}
	\Crefname{rem}{Remark}{Remarks}
	\Crefname{exmp}{Example}{Examples}
	\Crefname{notation}{Notation}{Notations}
	\Crefname{ass}{Assumption}{Assumptions}
	\usepackage[msc-links,initials]{amsrefs}
\makeatletter 
\def\@cite#1#2{[{\hbox{#1}\if@tempswa , #2\fi }]} 
\makeatother
	\newcommand{\lastindex}{\ell}

\newcommand{\N}{\mathbb{N}}

\newcommand{\Z}{\mathbb{Z}}
\newcommand{\Q}{\mathbb{Q}}

\newcommand{\F}{\mathbb{F}}

\newcommand{\Mat}[3]{\mathop{\mathrm{Mat}_{#1}^{#2}(#3)}}

\newcommand{\tuple}[1]{\mathbf{#1}}

	\newcommand{\rk}{\mathop{\mathrm{rk}}}

\newcommand{\cone}[2]{\mathop{W_{#1}(#2)}}

\newcommand{\ring}{R}
\newcommand{\ringint}{\mathcal{O}}
\newcommand{\primeideal}{\mathfrak{p}}
\newcommand{\qfield}{\mathbb{F}_q}
\newcommand{\numfield}{K}
\newcommand{\completion}{\mathfrak{o}}

        \newcommand{\alg}[1]{\mathfrak{g}_{#1}}
        
        \newcommand{\alginf}{\alg{}}
        \newcommand{\basis}{\mathcal{B}}
		\DeclareMathOperator{\alggroup}{\mathbf{G}}
		
			\newcommand{\GL}[2]{\mathop{\mathrm{GL}_{#1}(#2})}

\newcommand{\gp}[1]{\frac{#1}{1 - #1}}
		
\usepackage[final]{todonotes}

\newcommand{\tgroup}[1][]{$\mathcal{T}_{#1}$-group}

\usepackage{nicematrix}
\usepackage{booktabs}
\newcommand{\lattice}{\Lambda}
\newcommand{\locallie}{\mathfrak{g}}
\newcommand{\localcen}{\mathfrak{z}}
\renewcommand{\basis}{\mathbf{e}}
\newcommand{\basisder}{\mathbf{f}}
\newcommand{\inder}{r}
\newcommand{\ccn}{\mathrm{k}}

\DeclareMathOperator{\maj}{maj}
\newcommand{\len}{\ell}
\DeclareMathOperator{\nega}{neg}
\DeclareMathOperator{\des}{des}
\title[Zeta functions of unipotent group schemes of type $G$]{Univariate and bivariate zeta functions of unipotent group schemes of type $G$.}
\author{Michele Zordan}
\address{Department of Mathematics, Imperial College, London SW7 2BZ, UK}
\email{mzordan@imperial.ac.uk}
\begin{document}
\maketitle
\begin{abstract}
We compute the representation and class counting zeta functions for a family of torsion-free finitely generated nilpotent 
groups of nilpotency class $2$. These groups arise from a generalisation of one the families of unipotent groups schemes treated by 
Stasinski and Voll in \cite{stavol2011nilpotent,stavol2013hyper} 
and Lins in \cite{lin2019bivFGH}. 
The univariate zeta functions are obtained by specialising the respective bivariate zeta functions 
defined by Lins in \cite{lin2019arith}. These are also used to deduce a formula for a joint distribution on Weyl groups of type $B$.
\end{abstract}
\section{Introduction}
\thispagestyle{empty}
\subsection{A family of group schemes}
Stasinski and Voll \cite[Section~2.4]{stavol2011nilpotent} 
associate a unipotent group scheme $\alggroup_\Lambda$ with a $2$-nilpotent Lie lattice $\Lambda$ over the 
ring of integers of a number field. In this note we shall consider group schemes arising from the following $\Z$-Lie lattices. 
\begin{defn}
\label{def:groups_G}
Fix throughout  $m,n \in \N$. We define the following $\Z$-Lie lattice:
\begin{equation*}
\Lambda_{m, n}=\langle x_1,\dots, x_{m+n},z_{ij}, 1\leq i\leq m, 1\leq j\leq n\mid [x_i, x_{m+j}]=z_{ij}\rangle,
\end{equation*}
where all non-specified Lie brackets which do not follow formally from the given ones are zero.
\end{defn}
\subsection{Main results}
For ease of notation, we denote the group scheme $\alggroup_{\Lambda_{m,n}}$ by $G_{m, n}$. 
Without loss of generality we may assume that $m \leq n$. 
Indeed the $\Z$-Lie lattices $\Lambda_{m, n}$  and $\Lambda_{n,m}$ are isomorphic, so $G_{m, n}$ and 
$G_{n,m}$ are isomorphic as $\Z$-group schemes. Hence all formulas valid for $m \leq n$ work for $n \leq m$, 
after swapping $m$ and $n$.  Throughout, a \tgroup\ is torsion-free finitely generated and nilpotent.
\subsubsection{Twist representation zeta functions}
Let $G$ be a group and let $i \in \N$. We define $\mathrm{Irr}_i(G)$ to be the set of characters of the 
irreducible complex representations of $G$ of dimension $i$. If $G$ is a topological group, we stipulate that 
$\mathrm{Irr}_i(G)$ contains only the characters of the continuous irreducible representations.\par
Let $\theta_1, \theta_2 \in \mathrm{Irr}_i(G)$. We say that $\theta_1$ and $\theta_2$ are twist equivalent and 
write $\theta_1 \sim \theta_2$, when there is $\chi \in \mathrm{Irr}_1(G)$ such that $\theta_1 = \chi \theta_2$.
We define 
	\[
		\tilde{a}_i(G) = \# \left( \mathrm{Irr}_i(G)/ \sim \right).
	\] 
If the group $G$ is such that $\tilde{a}_i(G) \in \N$ for all $i \in \N$, we say that the group is {\em twist rigid}  and 
we define the twist representation zeta function of $G$ as	
	\[
		\zeta_G(s) = \sum_{i \in \N} \tilde{a}_i(G) i^{-s}.
	\]
It is well known that \tgroup s are twist rigid (cf.\ \cite[Theorem~6.6]{lubmag1985varieties}).
\begin{main}
\label{main:global}
Let $\ringint$ be the ring of integers of a number field $\numfield$, and let $\zeta_\numfield$ be the Dedekind zeta function
of $\numfield$. Then
	\[
		\zeta_{G_{m, n}(\ringint)}(s) = \prod_{i = 0}^{m - 1} \frac{\zeta_\numfield( s - n - i)}{\zeta_\numfield( s - i)}.
	\]
\end{main}
We briefly digress from the presentation of the main results to discuss the context of \Cref{main:global}. First of all, 
this theorem compares to \cite[Theorem B]{stavol2011nilpotent} where the authors compute the representation zeta function 
for three infinite families of \tgroup s generalising the Heisenberg group. The Lie lattice $\Lambda_{m, n}$
coincides with their $\mathcal{G}_{n}$ when $m = n$. Hence \cref{main:global} is a generalization of \cite[Theorem B]{stavol2011nilpotent}
for the groups  $G_n(\ringint)$. In particular the properties of the Dedekind zeta function imply that the results on functional 
equation and abscissa of convergence in \cite[Corollary 1.3]{stavol2011nilpotent} hold, mutatis mutandis, for the representation 
zeta function of $G_{m, n}(\ringint)$. Also the statement on meromorphic continuation holds but it is now a consequence of 
the more general \cite[Theorem~A]{dunvol2017uniform} by Dung and Voll. \par
Secondly, the groups $G_{m, n}(\ringint)$ are a generalization of the Grenham's groups $G_{1,n}$. These were used by 
Snocken to prove that every rational number may be attained as the abscissa of convergence of the representation zeta function 
of a \tgroup\ of class $2$ (cf.\ \cite[Theorem 4.22]{sno2012thesis}). In particular, \cref{main:global} may be used to
recover Snocken's result. Indeed, if $\times^k_Z G_{m, n}$ denotes the $k$-fold central product of $G_{m, n}$, then it is a 
well known fact that
	\[
		\zeta_{\times^k_Z G_{m, n}(\ringint)}(s)	= \zeta_{G_{m, n}(\ringint)}(k s).
	\]
This clearly has abscissa of convergence $(m + n) / k$. Notice that by \cite[Theorem~1.5]{hrumar2018definable} the 
abscissa of convergence must be rational (cf. also \cite[Theorem A]{dunvol2017uniform} for \tgroup s and \cite{zap2006zeta,stazor2018rational} 
for compact $p$-adic analytic groups and \cite{avn2011rationalabs} for arithmetic groups).\par
Finally, on the one hand, \Cref{main:global} also compares to \cite[Section~1.3]{carschvol2017traceless}, where Carnevale, Shecht\-er and Voll consider
Lie lattices obtained from $\mathcal{G}_n$ by adding an extra linear relation for its generators. On the other hand, in the 
same vein of generalising the Heisenberg group, \cref{main:global} (and its corollaries) compare with the result in 
\cite[Section~3.6]{vol2019delta} by Voll.
\subsubsection{Local representation zeta function}
Let $\numfield$ be a number field with ring of integers $\ringint$. For a non-zero prime ideal $\primeideal$ of $\ringint$ 
we denote  the completion of $\ringint$ at $\primeideal$ by $\ringint_\primeideal$.
Let $\Lambda$ be a $2$-nilpotent $\ringint$-Lie lattice. By \cite[Proposition~2.2]{stavol2011nilpotent} one has the following 
Euler factorization over the non-zero prime ideals of $\ringint$:
	\begin{equation}
	\label{eq:euler}
		\zeta_{\alggroup_\Lambda(\ringint)}(s)	=	\prod_{{\primeideal}} \zeta_{\alggroup_\Lambda(\ringint_{\primeideal})}(s).
	\end{equation}
Note that here we view $\alggroup_{\Lambda}(\ringint_\primeideal)$ as a topological group. Thus, by our convention, the factors 
on the right-hand side of the last equality ({\em local factors}) are defined by counting continuous representations only.
Note also that on the one hand, in its current form, 
equation \cref{eq:euler} applies to group schemes arising from Lie lattices over $\ringint$. On the other hand, 
$\alggroup_{\Lambda_{m , n} \otimes_{\Z} \ringint}(\ringint) \cong G_{m,n}(\ringint)$. To avoid 
unnecessary technicalities, now that $\ringint$ is fixed, we re-define $G_{m, n} = \alggroup_{\Lambda_{m , n} \otimes_{\Z} \ringint}$.
\begin{notation} 
Henceforth $X$, $Y$ and $Z$ will denote indeterminates in the field $\Q(X,Y,Z)$. We define the following objects.
	\begin{enumerate}
		\item Let $N\in \N$, we write
				\begin{align*}
				(\underline{N})_X 		&= (1 -  X^N)											&(\underline{0})_X  &= 1\\
				{(\underline{N})_X}!		&= (\underline{1})_X (\underline{2})_X \cdots (\underline{N})_X		&(\underline{0})_X! &= 1.
				\end{align*}
		\item	{\em Gauss polynomial}. For $a,b\in\N_0$ $b\leq a$, we define
				\[
					{a \choose b}_X = \frac{(\underline{a})_X!}{(\underline{a - b})_X! (\underline{b})_X!}.
				\]
		\item Let $j\in \N$. We write $[j]$ for the set $\lbrace 1, \dots, j\rbrace$ and $[j]_0$ for $\lbrace 0, \dots, j\rbrace$.
		\item We write $\lbrace i_1,\dots, i_\ell\rbrace_<$ for an {\em ordered subset} of $\N$, i.e.\ 
			a subset such that  $i_1< i_2 < \cdots < i_\lastindex$.
		\item {\em $X$-multinomial coefficient}. Let $j \in \N$ and let $I =\lbrace  i_1,\dots, i_\lastindex\rbrace_< \subseteq [j-1]_0$. 
			We write 
				\[
					{ j \choose I} = { j \choose i_\lastindex}_X {i_\lastindex \choose  i_{\lastindex -1}}_X \cdots {i_{2} \choose i_1}_X.
				\]
		\item {\em Pochhammer symbol.} Let $k\in \N$. We define
				\[
				(X ; Y)_k = \prod_{i = 0}^{k-1} (1 - XY^{i}).
				\]
		\end{enumerate}
\end{notation}
We introduce the following polynomial, which is related to counting $m \times n$ matrices of a given rank 
(see \cref{lem:rank} and \cref{eq:alt_form} for more detail).
\begin{defn}
Let $I=\lbrace i_1,\dots, i_\lastindex\rbrace_< \subseteq[m-1]_0$. We define
	\[
		f_{m, n}^I(X) = {m\choose I}_{X}(X^{n - m}X^{i_1 + 1};X)_{m - i_1}.
	\]
\end{defn}
\begin{main}
\label{main:additive_G}
Let $\primeideal$ be a non-zero prime ideal of $\ringint$ and $q$ be the cardinality of its residue field. Then
	\begin{equation}
		\label{eq:04}
		\zeta_{G_{m, n}(\ringint_\primeideal)}(s)=\sum_{I\subseteq[m-1]_0}f_{m, n}^I(q^{-1})\prod_{i\in I}\frac{q^{(m-i)(n+i)-
		(m-i)s}}{1-q^{(m-i)(n+i)-(m-i)s}},
	\end{equation}
where $I$ under the summation symbol denotes an ordered subset.
\end{main}
\subsubsection{Topological zeta functions} In \cite{ros2015topological} Rossmann defines topological representation zeta 
functions for \tgroup s arising from unipotent group schemes. 
\begin{cor}
\label{cor:topological}
For all non-zero prime ideals $\primeideal$ of $\ringint$, the topological representation zeta function of $G_{m, n}(\ringint_\primeideal)$ is 
	\[
		\prod_{i = 0}^{m - 1}\frac{s - i}{s - n - i}.
	\]
\end{cor}
\subsection{Bivariate zeta functions}
We shall prove \cref{main:global,main:additive_G} by computing the bivariate representation zeta function (cf.\ \cref{rem:BtoA}) 
introduced by Lins in \cite{lin2019arith}. In addition, we shall also compute another bivariate zeta function defined by Lins, 
which will give the class counting zeta function of $G_{m,n}(\ringint)$. 
We need to introduce some additional notation, in order to state the results of our computations.\par
Let $G$ be a group and let $n \in \N$. We define the following sequences:
	\begin{align*}
		r_n(G) 	&= \# \{ \text{irreducible complex characters of } G \text{ of degree } n\}\\
		c_n(G)	&=\# \{ \text{conjugacy classes of } G \text{ of size } n \}
	\end{align*}
Let $G$ be such 
that $r_n(G) \in \N$ or $c_n(G)\in \N$, for all $n \in \N$, we define respectively 
	\[
		\zeta_G^{irr}(s) = \sum_{n = 1}^{\infty} r_n(G) n^{-s}\qquad  \text{or}\qquad 
											\zeta_G^{cc}(s) = \sum_{n = 1}^{\infty} cc_n(G) n^{-s}.
	\]
Let $\alggroup$ be a unipotent group scheme over $\ringint$. Then, for any non-zero ideal $I \trianglelefteq \ringint$, $\alggroup(\ringint/I)$ 
is a finite group. Thus $r_n(\alggroup(\ringint/I))$ and $c_n(\alggroup(\ringint/I))$ are finite for all $n \in \N$.
\begin{defn}[{\cite[Definition~1.2]{lin2019bivFGH}}]
	\label{def:bivariate}
	The {\em bivariate representation zeta function} and the  {\em bivariate conjugacy class zeta function} of $\alggroup(\ringint)$
	are, respectively, 
		\begin{align*}
			\mathcal{Z}^{\mathrm{irr}}_{\alggroup(\ringint)}(s_1, s_2)	
					&	=	 \sum_{(0) \neq I \trianglelefteq \ringint} \zeta^{\mathrm{irr}}_{\alggroup(\ringint/I)}(s_1) \lvert \ringint : I \rvert^{-s_2}\\
			\mathcal{Z}^{\mathrm{cc}}_{\alggroup(\ringint)}(s_1, s_2)	
					&	=	 \sum_{(0) \neq I \trianglelefteq \ringint} \zeta^{\mathrm{cc}}_{\alggroup(\ringint/I)}(s_1) \lvert \ringint : I \rvert^{-s_2}.		
		\end{align*}
\end{defn}
By \cite[Proposition~2.4]{lin2019arith}, bivariate zeta functions have an Euler product factorisation. Namely, 
for $* \in \{ \mathrm{irr}, \mathrm{cc}\}$,
	\begin{equation}
		\label{eq:euler_biv}
		\mathcal{Z}_{\alggroup(\ringint)}^{*}(s_1, s_2) = \prod_{\primeideal \in \mathrm{Spec}(\ringint) \smallsetminus \{ (0)\}}
											\mathcal{Z}_{\alggroup(\ringint_\primeideal)}^{*}(s_1, s_2).
	\end{equation}
\newcommand{\zirr}{\mathcal{Z}_{G_{m,n}(\ringint_\primeideal)}^{\mathrm{irr}}(s_1, s_2)} 
\newcommand{\zirrcomp}{\mathcal{Z}_{G_{m,n}(\completion)}^{\mathrm{irr}}(s_1, s_2)} 
\newcommand{\zcc}{\mathcal{Z}_{G_{m,n}(\ringint_\primeideal)}^{\mathrm{cc}}(s_1, s_2)} 
\begin{main}
\label{main:bivariate}
Let $\primeideal$ be a non-zero prime ideal of $\ringint$ and $q$ be the cardinality of its residue field. 
Define
	\begin{align*}
		\delta_{n,m} 		&= n - m,	\\
		\overline{a}_{m,n}(i)	&= (m - i) ( n + i) + 2i + \delta_{n,m}.
	\end{align*}
Then
	\begin{align*}
		\zirr			&= \frac{1}{1 - q^{\overline{a}_{m,n}(m) - s_2}} \sum_{I \subseteq [m - 1]_0}
									\hspace{-2pt}															
									f_{m,n}^I(q^{-1}) \prod_{i \in I}
									\frac{q^{\overline{a}_{m,n}(i)  - (m - i) s_1 - s_2}}{1 - q^{\overline{a}_{m,n}(i)  - ( m - i) s_1 - s_2}}.		\\
		\zcc 			&= \frac{N^{\mathrm{cc}}_{m,n; q}(q^{-s_1}, q^{-s_2})}{D^{\mathrm{cc}}_{m,n; q}(q^{-s_1}, q^{-s_2})},
	\end{align*}
where
	\begin{align*}
N^{\mathrm{cc}}_{m,n; q}(T_1, T_2)		=	&\,T_{1}q^{m + n}  \\
						&+ ( T_{1}^{2m + n} T_{2}^{2}  \, q^{mn + 1} - T_{1}^{m + 1} T_{2}) \,(q^{m} + q^{n}) \,q^{mn}\\
						&+ ( T_{1}^{m + n + 1} T_{2}^{2}\, q^{mn} - T_{1}^{m + n}  T_{2}\, q)\, (q^{m} + q^{n} -1)\,q^{mn}\\
						&-T_{1}^ {2( m + n)}T_{2}^{3} \, q^{3 mn + 1}\\
%
D^{\mathrm{cc}}_{m,n; q}(T_1, T_2)		=\,	&T_1q^{m + n} (1 - T_1^{m + n - 1} T_2 q ^{mn + 1}) (1 - T_1^{n} T_2  q^{(m - 1)(n + 1) + 1}) \\
						&\cdot (1 - T_1^{m} T_2  q^{(n - 1)(m + 1) + 1}) ( 1 - T_2 q^{mn}).	
	\end{align*}
\end{main}
The second part of this theorem is proved in \cite[Section~3.3]{lin2019bivFGH} when $m = n$. The proof we shall 
give in \Cref{sec:cc}, however, is not a direct generalisation of that proof.
\subsubsection{Class counting zeta function}
As anticipated, we shall use \cref{main:bivariate} to deduce the main results on the representation zeta function from its 
bivariate counterpart. The bivariate conjugacy class zeta function will, instead, be used to obtain 
a special case of \cite[Proposition~8.7]{rosvol2019graphs} (for $r=2, \allowbreak n_1=n, n_2=m$).\par
Recall that the number of conjugacy classes of a group $G$ is called the {\em class number} and is denoted by $\ccn(G)$. 
Let $\alggroup$ be a unipotent group scheme over $\ringint$. The class counting zeta function of $\alggroup(\ringint)$ is defined as
	\[
		\zeta_{\alggroup(\ringint)}^{\ccn}(s) = \sum_{(0) \neq I \trianglelefteq \ringint} \ccn(\alggroup(\ringint/I)) \lvert \ringint : I \rvert^{-s}.
	\]
Unipotent groups have the strong approximation property (cf.\ \cite[Lemma~5.5]{plarap1994alggroup}) 
and therefore, by \cite[Lemma~8.1]{berderonnpir2013uniform},
the class counting zeta function has an Euler product decomposition.
\begin{cor}
\label{main:class_number}
The class counting zeta function $\zeta_{G_{m,n}(\ringint)}^{\ccn}(s)$ is
	\[
		\prod_{\primeideal}
		\frac{q^{2  m n -m - n + 1- 2s} - ((q + 1) q^{m} - (q^{m} - q - 1) q^{n} - q) q^{m n - m - n - s} + 1}
		{(1 - q^{m n + m - n - s} ) (1 - q^{m n + n - m - s} ) (1 - q^{m n + 1 - s} )},
	\]
	where $q$ is the residue field cardinality of $\primeideal$, and $\primeideal$ ranges over the non-zero prime ideals of $\ringint$.
\end{cor}
\begin{rem}
Although we shall use the bivariate conjugacy class zeta function to derive the class counting zeta function, it is worth 
noting here, that the latter is also part of the larger theory of {\em ask zeta functions} formulated by Rossmann in \cite{ros2018average_i} (see Theorem~1.17 therein). 
\end{rem}
\subsubsection{Statistics on Weyl groups}
The last main result in this work follows from an alternative expression of the bivariate representation zeta function in Theorem~\ref{main:bivariate}. 
We recall some definitions and basics on Weyl groups from \cite[Section~4.2]{stavol2011nilpotent}. 
Let $(W,S)$ be a finite Coxeter system, with Coxeter group $W$ and Coxeter generating set $S$.
\begin{notation}
Let $w\in W$. We denote the Coxeter length of $w$ by $\len(w)$. That is, $\len(w)$ is the shortest length of a word in the elements of 
$S$ representing $w$. We denote the (right) descent type of $w$ by 
	\[
		D(w) = \{  s \in S \mid \len(ws) < \len(w) \}.
	\]
\end{notation}
We have some additional notation for Weyl groups of type $B$. Let $[\pm m]_0 = \{ -m, \dots, 0, \dots, m\}$. 
The group $B_m$ is defined as the subgroup of the symmetric group on $[ \pm m]_0$, consisting of those permutations $w$ such that, 
for all $i \in [ \pm m]_0$, $w(- i ) = - w(i)$.  We choose the standard set of Coxeter generators for $B_m$, namely 
$S = \{ \sigma_0, \sigma_1, \dots, \sigma_{n - 1} \}$ where
	\begin{align*}
		s_0 	&= (-1, 1)\\
		s_i	&= (-i - 1, -i) (i, i + 1)  \text{ for all } i \in [m -1].
	\end{align*}
For $w \in B_m$, we define
	\[
		\nega(w) = \# \{ i \in [m] \mid w(i) < 0 \}.
	\]
Applying \cite[Lemmas~4.4 and~4.5]{stavol2011nilpotent} to the formula in Theorem~\ref{main:additive_G}, 
we obtain the following expression of the bivariate representation zeta function in the fashion of \cite[Lemma~5.4]{lin2019bivFGH}.
\begin{cor}
\label{cor:weyl}
	Let $\primeideal$ be a non-zero prime ideal $\primeideal$ in $\ringint$ and let $q$ be its residue field cardinality. 
	For $w \in B_m$, let
		\begin{align*}
			\chi_{m,n}(w) 	&=  (-1)^{\nega(w)}\\
			h_{m,n}(w)	&= \len(w) + \delta_{n,m} \nega(w).
		\end{align*}
	Then
		\[
			\mathcal{Z}_{G_{m,n}(\ringint_\primeideal)}^{\mathrm{irr}}(s_1, s_2) = \frac{	\sum_{w \in B_m} \chi_{m,n}(w)  q^{-h_{m,n}(w)} \prod_{i \in D(w)} 
																	q^{\overline{a}_{m,n}(i)  - ( m - i) s_1 - s_2}}
										{\prod_{i = 0}^{m} ( 1 - q^{\overline{a}_{m,n}(i)  - ( m - i) s_1 - s_2} ) }.
		\]
\end{cor}
Note that -- in comparison to \cite[Lemma~4.4]{stavol2011nilpotent} -- the formula in the last corollary has an extra factor for $i = m$ in the denominator. This accounts for the factor $(1 - q^{\overline{a}(m) - s_2})$ in the expression for 
$\zirr$ in \cref{main:bivariate}.\par
Similar to \cite[Propositions~5.5 and 5.6]{lin2019bivFGH} we deduce the following formula for a joint distribution of the following 
statistics on $B_m$. We define

	\begin{align*}
		\sigma_{mn} 	& = \sum_{i \in D(w)} (m - i) (n + i)\\
		\maj(w)		& = \sum_{i \in D(w)} i\\
		\des(w) 		& = \lvert D(w) \rvert.
	\end{align*}
\begin{cor}
\label{cor:joint}
For all $m, n \in \N$ with $m \leq n$, we have the following identity in $\Q[X, Z]$. 
	\begin{multline*}
		\sum_{w \in B_m} (-1)^{\nega(w)} X^{-(\sigma_{mn} - \ell + 2 \maj)(w) + \delta_{mn}(\des - \nega) (w)} Z =\\
					\prod_{i = 3}^{n}(1 - X^{(m - i) (n + i) + 2i + \delta_{mn}} Z) \cdot \\
		(X^{2mn - m - n + 1} Z^2 - ((X + 1) X^m - (X^m - X - 1) X^n - X)) X^{mn -m - n} + 1).
	\end{multline*}
\end{cor}	

\begin{notation}
Throughout, $\Z$ denotes the integers, $\N$ the set of positive integers and $\N_{0} = \lbrace 0\rbrace \cup \N$ the set 
of natural numbers. The set of rational numbers is denoted by $\Q$.\par
Let $\ring$ be a ring, an $\ring$-Lie lattice is a free finitely generated $\ring$-module endowed with a Lie bracket. 
If $\alginf$ is an $\ring$-Lie lattice, we write $\alginf'$ for its derived Lie sublattice. If not otherwise specified, when 
$\ring$ is considered as an $\ring$-Lie lattice it is always endowed with the trivial Lie bracket.\par
Let $i,j \in\N$. The ring of $i \times j$ matrices with entries in $\ring$ is denoted by $\Mat{i \times j}{}{\ring}$. The zero 
element of this ring is denoted by $O_{i \times j}$. If $i = j$
the identity matrix is denoted by $\mathrm{\mathrm{Id}}_{j}$. The diagonal matrix with $a_1,\dots,a_j \in \ring$ on the diagonal 
is denoted by $\mathrm{diag}(a_1,\dots,a_j)$.\par
We shall denote tuples by $\tuple{x}, \tuple{y}, \tuple{z}, \dots$, while their components will be
denoted by $x_1, x_2, x_3,\dots,y_1, y_2, y_3,\dots,z_1, z_2, z_3,\dots$ respectively. Often we shall represent 
$mn$-tuples of elements in a ring $\ring$ as $m \times n$ matrices. In this case the components 
of the $mn$-tuple $\tuple{x}\in\Mat{m \times n}{}{\ring}$ are denoted by $x_{ij}$ for $i \in \lbrace 1,\dots, m\rbrace$
and $j \in \lbrace 1,\dots, n\rbrace$.
\end{notation}
\subsection{Acknowledgement}
The formulae for the (univariate) representation zeta functions computed here were first obtained during 
my doctoral studies. I wish to thank my PhD supervisor Christopher Voll and I acknowledge 
financial support by the School of Mathematics of the University of Southampton, the Faculty of Mathematics of the University 
of Bielefeld, and CRC 701. The remaining parts of this work were completed while supported by Research Project G.0939.13N and G.0792.18N of the 
Research Foundation - Flanders (FWO), the University of Auckland, and Imperial College London. I am grateful to Ben Martin  and 
Christopher Voll for their comments.
\section{Preliminaries on bivariate zeta functions}
Let $\primeideal$ be a non-zero prime ideal in $\ringint$. Lins shows in \cite{lin2019bivFGH} that the 
$\primeideal$-local factor in \eqref{eq:euler_biv} may be expressed in terms of $\primeideal$-adic integrals. 
We briefly recall how it is done. 
\subsubsection{Igusa-type integrals}
For the rest of this note, let $\completion = \ringint_\primeideal$. 
The main tool is an Igusa-type integral associated with a matrix of linear forms. Namely, 
let $d \in \N$ and let $\tuple{X} = (X_1, \dots, X_d)$. Let $\mathcal{R}(\tuple{X})$ be a matrix of linear forms in 
$\completion[\tuple{X}]$. Let $\numfield_\primeideal$ be the field of fractions of $\completion$. We define
	\begin{equation}
	\label{eq:u_R}
		u_{\mathcal{R}}=  {\rk}_{\numfield_\primeideal(\tuple{X})} (\mathcal{R}(\tuple{X})).
	\end{equation}
Let $\mu$ be the additive Haar measure on $\completion^{d + 1}$. For $k \in [u]$, let also 
	\[
		\mathcal{F}^{k}(\mathcal{R}(\tuple{X}))
	\] 
be the ideal of $\completion[\tuple{X}]$ generated by the $k$-minors of $\mathcal{R}(\tuple{X})$. 
We set
	\[
		\cone{d}{\completion} 	= \completion^d \smallsetminus \primeideal^{d}.	
	\]
One defines
	\[
		\mathcal{Z}_{\mathcal{R}}(\rho, \tau) = \frac{1}{1 - q^{-1}} \int_{(w, \tuple{x}) \in \primeideal \times \cone{d}{\completion}}
			\lvert w \rvert^{\tau}_\primeideal 
				\prod_{k =1}^{u_{\mathcal{R}}} \frac{\|\mathcal{F}^k (\mathcal{R}(\tuple{x}))\cup w\mathcal{F}^{k-1}(\mathcal{R}(\tuple{x})) \|_\primeideal^{\rho}}
				{\| \mathcal{F}^{k - 1} (\mathcal{R}(\tuple{x}))\|_\primeideal^{\rho}}\, d\mu.
	\]
\subsubsection{Commutator matrices}

The matrices of linear forms one uses in this context are related to the structure constants of the Lie lattice defining the unipotent 
group scheme. Recall that $\lattice$ is a Lie lattice over $\ringint$. Let
	\begin{align*}
		\locallie	&= \lattice \otimes_{\ringint} \completion		&	\locallie' 	& = [\locallie, \locallie]	& \localcen = Z(\locallie),
	\end{align*}
	\begin{align*}
		   	h	&= {\rk}_{\completion} (\locallie)	,		& a		&= {\rk}_{\completion} (\locallie/\localcen),		
		   &	b	&= {\rk}_{\completion} (\locallie'),		& \inder	&= {\rk}_{\completion} (\locallie/\locallie'),		
		   &	z	&= {\rk}_{\completion} (\localcen).
	\end{align*}
Let $\basis = (e_1 + \localcen, \dots, e_a + \localcen)$ be an ordered set of $\completion$-module generators of $\locallie/ \localcen$ and 
$\basisder =  (f_1, \dots, f_b)$ be an ordered set of $\completion$-module generators of $\locallie'$.
For $i, j \in [a]$, the following equation defines the {\em structure constants} of $\locallie$: $\lambda_{ij}^1, \dots, \lambda_{ij}^b\in \completion$ such that 
	\[
		[e_i, e_j] = \sum_{k = 1}^{b}  \lambda_{ij}^k f_k.
	\]
	
\begin{defn}[{\cite[Definition~2.1]{brivol2015enumerating}}]
\label{def:cmatrix}
Let $\tuple{X} = (X_1, \dots, X_a)$ and $\tuple{Y} = (Y_1, \dots, Y_b)$ be independent variables. We define the following
{\em commutator matrices} of $\completion$-linear forms:
(with respect to $\basis$ and $\basisder$) 
	\begin{align*}
		A(\tuple{X}) 	&\in \Mat{a\times b}{}{\completion[\tuple{X}]} 	&	\text{ where }	A(\tuple{X})_{ik} &= \sum_{j = 1}^a \lambda_{ij}^k X_j , 
	&	i\in [a], k \in [b],	\\
		B(\tuple{Y}) 	&\in \Mat{a\times a}{}{\completion[\tuple{Y}]} 	&	\text{ where }	B(\tuple{Y})_{ij}  &= \sum_{k = 1}^b \lambda_{ij}^k Y_k  , 
	&	i, j\in [a].	\\
	\end{align*}
\end{defn}
We used the structure constants to define the two commutator matrices. However, the 
duality between them may be expressed in a coordinate-free way. This has been done by Rossmann in \cite[Remark~4.13]{ros2020average_ii}.
\subsubsection{Bivariate zeta functions as integrals}
The following result by Lins gives the bivariate zeta functions in terms of $\primeideal$-adic integrals. Note that here $u_B$ is defined as 
in \eqref{eq:u_R} and is twice $u_{B_\Lambda}$ from  \cite{lin2019arith}.
\begin{prop}[{\cite[Proposition~4.8]{lin2019arith}}]
\label{prop:padic}
We may express the bivariate zeta functions of $\alggroup_\lattice(\completion)$ as
	\begin{align*}
		\mathcal{Z}^{\mathrm{irr}}_{\alggroup_\lattice(\completion)}(s_1, s_2)	
				&	=	\frac{1}{1 - q^{\inder - s_2}}	 
						\Bigg( 1 + \mathcal{Z}_{B}\bigg( -\frac{s_1 + 2}{2}, \frac{s_1 + 2}{2} \, u_B + s_2 - h - 1 \bigg) \Bigg)\\
		\mathcal{Z}^{\mathrm{cc}}_{\alggroup_\lattice(\completion)}(s_1, s_2)	
				&	=	 \frac{1}{1 - q^{z - s_2}}	 
						\bigg( 1 + \mathcal{Z}_{A}\Big( -(s_1 + 1), (s_1 + 1) u_A + s_2 - h - 1 \Big) \bigg).	
	\end{align*}
\end{prop}
\begin{rem}
In \cite{lin2019arith}, the last proposition is stated with restrictions on the residue field characteristic in relation to the nilpotency class 
of $\lambda$. However, in our case -- and whenever the nilpotency class is $2$ -- it is possible to use the construction of $\alggroup_\Lambda$ 
introduced in \cite[Section~2.4]{stavol2011nilpotent}. The Kirillov orbit method explained there gives the first equality 
in the statement of Proposition~\ref{prop:padic} for general $p$. Moreover, for all $x, y \in \alggroup_\Lambda(\completion)$, 
	\[
		xyx^{-1}y^{-1} = [x , y].
	\] 
Thus $C_{\alggroup_\Lambda(\completion)} (x) = \ker (\mathrm{ad}_x$). Following the proof of 
Proposition~4.8 in \cite{lin2019arith}, this gives the second equality for general $p$. 
Indeed, the result above appears without restriction on the prime 
for nilpotency class $2$ as \cite[Proposition~2.6]{lin2019bivFGH}.
\end{rem}
\section{Computing the bivariate representation zeta function}
In this section we compute the $\primeideal$-adic  integral giving the local bivariate representation zeta function of 
$G_{m,n}(\completion)$. 
For convenience, in what follows we identify $mn$-tuples with  
$m \times n$ matrices, so that, for example, $\completion^{mn}$ is identified with $\Mat{m\times n}{}{\completion}$.\par
We start by computing the matrix $B$ for the Lie lattice $\lattice_{m,n}$. Using 
$(x_1+ \mathfrak{z}, \dots, x_{m + n} + \mathfrak{z})$ as $\basis$ and $(z_{11}, z_{12}\dots, z_{mn})$ as $\basisder$, we 
immediately see that the $B$-commutator matrix with respect to $\basis$ and $\basisder$ is 
	\[
		B(Y_{11}, Y_{12}, \dots, Y_{m n}) = 	\begin{bNiceArray}{ccc|ccc}[margin]
				\Block{3-3}{O_{m \times n}}	&		&		& Y_{11}					&\Cdots	& Y_{1 n}\\
										&		&		&\Vdots					&		&\Vdots\\
										&		&		& Y_{m1}					&\Cdots	& Y_{m n}\\
				\hline			-Y_{11}	&\Cdots	& -Y_{m1}	&\Block{3-3}{O_{n \times m}}	&		&\\
								\Vdots	&		& \Vdots	&						&		&\\
								-Y_{1n}	&\Cdots	& -Y_{mn}	&						&		&
									\end{bNiceArray}.
	\]
It follows that $u_B = 2 m$ because $m \leq n$.\par
Let $N \in \N_0$  and let $\completion_N = \completion/ \primeideal^N$ (so that $\completion_0 = \completion/\completion$). 
By abuse of notation, in what follows
$B(\tuple{w})$ is the reduction of $B(\tuple{Y})$ modulo $\primeideal^N$ evaluated at 
$\tuple{w}\in(\completion/\primeideal^N)^{m n}$. 
\begin{defn}
For $N > 0$, we define
	\begin{align*}
		\nu_{N}(B(\tuple{w}))			&=(\min \lbrace a_i, N \rbrace)_{i \in \lbrace 1,\dots, m\rbrace},
	\end{align*}
where $\primeideal^{a_1}, \primeideal^{a_1}, \dots, \primeideal^{a_m}, \primeideal^{a_m}$ are the first 
$2m$ elementary divisors of a lift of $B(\tuple{w})$ to a matrix with entries in $\completion$. Note that $u_B = 2m$, so  $B(\tuple{w})$ has at most $2m$ non-maximal elementary divisors. We extend this definition 
for $N = 0$ by saying that $\nu_{0}(B(\tuple{w})) =
(0, \dots, 0) \in \N_0^m$, where $\tuple{w}$ is the unique element of $(\completion/\completion)^{mn}$.
\end{defn}
Let $I=\lbrace i_1,\dots,i_\lastindex \rbrace_{<}\subseteq [m - 1]_0$. We set $i_0=0$ and $i_{\lastindex + 1} = m$. Let
	\begin{align*}
		\mu_j		&=i_{j+1}-i_j 							&&\text{for } j\in \lbrace 0,\dots, \lastindex \rbrace, \\
		N			&=\sum_{i \in I} r_{i}				&&\text{for } \mathbf{r}_I=(r_{i_1},\dots,r_{i_\lastindex})\in \N^{ I}.
	\end{align*}
Let 
$\cone{m n}{\completion_N} = (\completion_N)^{m n} \smallsetminus (\primeideal/\primeideal^N)^{m n}$ for $N > 0$, and 
$\cone{m n}{\completion_0} = (\completion/ \completion)^{mn}$.
We define
	\begin{multline*}
		\mathrm{N}^{\completion}_{I,\mathbf{r}_I}(m,n) = \lbrace \tuple{w}\in \cone{m n}{\completion_N}\mid \nu_{N}(B(\tuple{w}))=(\underbrace{0,\dots,0}_{\mu_\lastindex},\underbrace{r_{i_{\lastindex}},\dots,r_{i_{\lastindex}}}_{\mu_{\lastindex-1}},\\
			\underbrace{r_{i_{\lastindex}} + r_{i_{\lastindex-1}} ,\dots,r_{i_{\lastindex}} + r_{i_{\lastindex-1}} }_{\mu_{\lastindex-2}},%
					\dots%
						,\underbrace{N,\dots,N}_{\mu_0})\in \N_0^m\rbrace.
	\end{multline*}
Note for later that we use the convention $\N^{\emptyset} = \{ \emptyset\}$. Hence, for $I = \emptyset$, we have 
$\mathbf{r}_I = \emptyset$. Thus,  in this case, $N = 0$ and 
$\lvert \mathrm{N}^{\completion}_{I,\mathbf{r}_I}(m,n) \rvert = 1$. \par
Following the proof of \cite[Equation~4.2]{lin2019bivFGH}, we obtain that
	\begin{multline}
	\label{eq:p_series}
		\mathcal{Z}_{G_{m,n}(\completion)}^{\mathrm{irr}}(s_1, s_2) = \\
		\frac{1}{1 - q^{m + n - s_2}}
		\sum_{I\subseteq [m - 1]_0}
			\sum_{\mathbf{r}_I\in\N^{I}} 
				\vert \mathrm{N}^{\completion}_{I,\mathbf{r}_I}(m,n)\vert
				\,q^{\sum_{i \in I} r_i (- ( m - i) s_1 - s_2 + 2i + \delta_{m,n})}.
	\end{multline}
\subsection{Lifting matrices of a given rank}
In what follows we compute the quantities $\vert \mathrm{N}^{\completion}_{I,\mathbf{r}_I}(m,n)\vert$.  The techniques are
similar to those in \cite{stavol2011nilpotent,zor2016adjoint}, in that they are based on lifting the Smith normal forms of  $B(\tuple{Y})$ evaluated over $\qfield$.\par 
Let $\pi$ be a generator of $\primeideal$. According to our convention of representing $mn$-tuples as matrices, for all $I \subseteq [m - 1]_0$ and $\mathbf{r}_I \in \N^I$, we identify the set $\mathrm{N}^{\completion}_{I,\mathbf{r}_I}(m,n)$ with the set of 
matrices over $\completion_N$ having Smith normal form 
	\[
		\mathrm{diag}(\underbrace{1, \dots, 1}_{\mu_\ell},\underbrace{\pi^{r_{i_1}}, \dots, \pi^{r_{i_1}}}_{\mu_{\ell - 1}}, \dots,
																	\underbrace{\pi^N, \dots, \pi^N}_{\mu_0}).
	\]

\subsubsection{Base step}
The computation is an induction on the size of the set $I$. We start by considering the base case $\lvert I \rvert = 1$. We need the following notation.
\begin{defn}
Let $N \in \N$ and $i\in[m]_0$. 
We define $ \Mat{m\times n}{m-i}{\completion_N}$ as the set of matrices in $\Mat{m\times n}{}{\completion_N}$ that 
have Smith normal form 
	\[
		\begin{bmatrix}
			\mathrm{Id}_{m - i}	&	0\\
			0		&	0
		\end{bmatrix}.
	\]
Note that, with the notation above, 
\begin{equation*}
\Mat{m\times n}{m-i}{\F_q}=\lbrace x\in \Mat{m\times n}{}{\F_q}\mid \rk(x)=m-i\rbrace. 
\end{equation*}
\end{defn}
We start by considering the situation over the residue field. The following lemma appears, with different notation, as Proposition~3.1 in 
\cite{laktho1994counting}.
\begin{lem}
\label{lem:rank}
For $i\in [m]_0$, 
	\begin{equation*}
	\vert\Mat{m\times n}{m-i}{\F_q}\vert={n\choose n-m+i}_{q^{-1}}(q^{-i-1};q^{-1})_{m-i}\cdot q^{(m-i)(n+i)}.
	\end{equation*}
\end{lem}
We now consider how matrices of $\Mat{m\times n}{m-i}{\F_q}$ lift to $\completion_N$ for $N \in \N$. The following lemma will form the base of our induction.
\begin{lem}
\label{lem:lift_rank}
Let $i\in[m]_0$ and $\overline{\tuple{x}} \in \Mat{m\times n}{m-i}{\F_q}$. Let also $N \in \N$. Then there are exactly 	
	\[
		q^{(N - 1) (m - i)(n + i)}
	\]
different lifts of $\overline{\tuple{x}}$ to $ \Mat{m\times n}{m-i}{\completion_N}$.
\end{lem}
	\begin{proof}
		The matrix $\overline{\tuple{x}}$ has $m - i$ independent rows, say $\overline{v}_1, \dots, \overline{v}_{m - i}$. Each of 
		these lifts in $q^{(N - 1) n}$ ways.\par
		Assume, now, that we have chosen $v_1, \dots, v_{m - i} \in (\completion_N)^n$ lifting  
		$\overline{v}_1, \dots, \overline{v}_{m - i}$ respectively. 
		The remaining $i$ rows of $\overline{\tuple{x}}$ are $\qfield$-linear combinations of the initial rows, and -- in order to achieve 
		the desired Smith normal form -- each of these linear combinations needs to be lifted to an $\completion_N$-linear 
		combination of $v_1, \dots, v_{m - i}$.\par
		Let $\overline{\alpha}_1, \dots \overline{\alpha}_{m - i} \in \qfield$ and let 
		$\overline{v} = \sum_{j = 1}^{m - i} \overline{\alpha}_j \overline{v}_j$. Since $\overline{v}_1, \dots \overline{v}_{m - i}$ are 
		linearly independent over $\qfield$, their lifts $v_1, \dots v_{m - i}$ are linearly independent over $\completion_N$. This 
		means that an $\completion_N$-linear combination of $v_1, \dots, v_{m - i}$, say $ \sum_{j = 1}^{m - i} \alpha_j v_j$, 
		is a lift of $\overline{v}$ if and only if $\alpha_j \equiv \overline{\alpha}_j\bmod \primeideal$ for all $j \in [m - i]$. There are exactly $q^{N - 1}$ 
		independent choices for each of the $\alpha_1, \dots, \alpha_{m - i}$. Thus a linear combination of 
		$\overline{v}_1, \dots, \overline{v}_{m - i}$ lifts to a linear combination of $v_1, \dots, v_{m - i}$ in $q^{(N - 1) (m - i)}$ different ways. 
		We conclude that there are 
			\[
				q^{(N - 1) (m - i) n + (N - 1) (m - i ) i} = q^{(N - 1) (m - i)(n + i)}
			\]
		lifts of $\overline{\tuple{x}}$ having the prescribed Smith normal form.
	\end{proof}
\subsubsection{Inductive step} The next lemma will be key for the inductive step of the proof. 
We need the following definition for ease of notation.
\begin{defn}
Let $I = \{i_1, \dots, i_\ell \}_{<} \subseteq [m -1]_0$, $I \neq \emptyset$. Let $\mathbf{r}_I \in \N^{I}$. We define
	\begin{align*}
		I^{-}				&= I \smallsetminus \{ i_\ell \},		\\
		\mathbf{r}_I^{-}		& = 	\begin{cases}
								\emptyset					&\lvert I \rvert = 1\\	
								(r_{i_1}, \dots, r_{i_{\ell - 1}}) 	& \lvert I \rvert > 1.
							\end{cases}
	\end{align*}
\end{defn}
\begin{lem}
\label{lem:ind_step}
Let $I = \{i_1, \dots, i_\ell \}_{<} \subseteq [m - 1]_0$ and $\mathbf{r}_I \in \N^{I}$. Then 
	\begin{equation*}
		\lvert \mathrm{N}^{\completion}_{I, \mathbf{r}_I}(m,n) \rvert =  \lvert \Mat{m\times n}{m - i_\ell}{\F_q}\vert \cdot 
		 \lvert \mathrm{N}^{\completion}_{I^{-}, \mathbf{r}_I^{-}}(i_\ell, n - m + i_\ell) \rvert 
		 \cdot q^{(N - 1) (m - i_{\ell}) (n + i_{\ell})}. 
	\end{equation*}
\end{lem}
	\begin{proof}
		Let $\overline{\tuple{x}} \in \Mat{m\times n}{m - i_\ell}{\F_q}$. In order to simplify the notation, in this proof we set
			\[
				\mathrm{N}_{\overline{\tuple{x}}} = \{ \tuple{w} \in  \mathrm{N}^{\completion}_{I, \mathbf{r}_I}(m,n)  
					\mid \tuple{w} \text{ reduces  to } \overline{\tuple{x}}\, \bmod \primeideal \}.
			\]
		We shall show that there is an onto map 
			\[
				\mathrm{N}_{\overline{\tuple{x}}} \longrightarrow \mathrm{N}^{\completion}_{I^{-}, \mathbf{r}_I^{-}}(i_\ell, n - m + i_\ell)
			\]
		that has fibres of cardinality $q^{(N - 1) (m - i_{\ell}) (n + i_{\ell})}$.\par
		Since $\overline{\tuple{x}}$ has $\qfield$-rank $m - i_\ell$, there are two matrices $\overline{P} \in \GL{m}{\qfield}$ and 
		$\overline{Q} \in \GL{n}{\qfield}$ such that 
			\[
				\overline{P} \overline{\tuple{x}} \overline{Q} = 	\begin{bmatrix}
														\mathrm{Id}_{m - i_\ell}	&	0\\
														0			&	0
													\end{bmatrix}.
			\]
		Let $P \in \GL{m}{\completion_N}$ and $Q \in \GL{n}{\completion_N}$ be lifts of $\overline{P}$ and $\overline{Q}$ respectively. The 
		sets $\mathrm{N}_{\overline{\tuple{x}}}$ and $P\mathrm{N}_{\overline{\tuple{x}}}Q$ have the same cardinality. 
		Moreover, if ${\mathbf{r}_{I}^{\pm}} =  \mathbf{r}_{I}^{-} + (r_{i_\ell}, 0, \dots, 0)$, then the sets 
		$\mathrm{N}^{\completion}_{I^{-}, \mathbf{r}_I^{-}}(i_\ell, n - m + i_\ell)$ and 
			\[
				\mathrm{N}^{-} := \pi^{r_{i_\ell}} \cdot \mathrm{N}^{\completion}_{I^{-}, \mathbf{r}_I^{\pm}}(i_\ell, n - m + i_\ell)
			\]
		have the same cardinality. We shall find an onto map
			\[
				\sigma_{\overline{\tuple{x}}} :  P\mathrm{N}_{\overline{\tuple{x}}}Q \longrightarrow \mathrm{N}^{-}
			\]
		and prove that its fibres have cardinality $q^{(N - 1) (m - i_{\ell}) (n + i_{\ell})}$.\par
		Before proceeding with the definition of $\sigma_{\overline{\tuple{x}}}$ we recall the following fact 
		about block matrices with an invertible diagonal block. Let 
			\begin{align*}
				A 	&\in \GL{m - i_\lastindex}{\completion_N}{},	& B	&\in \Mat{(m - i_\lastindex) \times (n - m + i_\lastindex)}{}{\completion_N},\\
				C	&\in \Mat{i_\lastindex \times (m - i_\lastindex)}{}{\completion_N},				&D	&\in \Mat{i_\lastindex \times (n - m + i_\lastindex)}{}{\completion_N}.
			\end{align*}
		We have that
			\begin{multline}
			\label{eq:block_mult}
				\begin{bNiceMatrix}
				\mathrm{Id}_{m - i_\lastindex}	&	O_{(m - i_\lastindex) \times  i_\lastindex}\\
				-CA^{-1}					&	\mathrm{Id}_{i_\lastindex}
				\end{bNiceMatrix}
				\begin{bNiceMatrix}
				A		& B\\
				C		& D
				\end{bNiceMatrix}
				\begin{bNiceMatrix}
				\mathrm{Id}_{m - i_\lastindex}					& -A^{-1}B \\
				O_{(n - m + i_\lastindex) \times  (m - i_\lastindex)}	& \mathrm{Id}_{n - m + i_\lastindex}
				\end{bNiceMatrix}															=\\
				\begin{bNiceMatrix}
				\mathrm{Id}_{m - i_\lastindex}	&	O_{(m - i_\lastindex) \times  i_\lastindex}\\
				-CA^{-1}					&	\mathrm{Id}_{i_\lastindex}
				\end{bNiceMatrix}
				\begin{bNiceMatrix}
				A		& O_{(m - i_\lastindex) \times (n - m + i_\lastindex)}\\
				C		& D - CA^{-1}B
				\end{bNiceMatrix}
				=\\
							\begin{bNiceMatrix}
				A									& O_{(m - i_\lastindex) \times (n - m + i_\lastindex)}\\
				O_{i_\lastindex \times (m - i_\lastindex)}		& D - CA^{-1}B
				\end{bNiceMatrix}.
			\end{multline}
		We shall now define $\sigma_{\overline{\tuple{x}}}$. Let $\tuple{w} \in P\mathrm{N}_{\overline{\tuple{x}}}Q$. Then 
		there are (uniquely determined)
			\[
				\begin{aligned}
					\tuple{a} 	&\in \GL{m - i_\lastindex}{\completion_N}{},					& \tuple{b}	&\in \Mat{(m - i_\lastindex) \times (n - m + i_\lastindex)}{}{\completion_N},\\
					\tuple{c}	&\in \Mat{i_\lastindex \times (m - i_\lastindex)}{}{\completion_N},	&\tuple{d}	&\in \Mat{i_\lastindex \times (n - m + i_\lastindex)}{}{\completion_N},
				\end{aligned}	
			\text{ such that }
				\begin{aligned}
					\tuple{w} &= 	\begin{bNiceMatrix}
								\tuple{a}	&	\tuple{b}\\
								\tuple{c}	&	\tuple{d}
							\end{bNiceMatrix}.
				\end{aligned} 		
			\]
		We set $\sigma_{\overline{\tuple{x}}}(\tuple{w}) = \tuple{d} - \tuple{c}\tuple{a}^{-1} \tuple{b}$. This gives a function 
		$P\mathrm{N}_{\overline{\tuple{x}}}Q \rightarrow \mathrm{N}^{-}$, because 
		$\tuple{d} - \tuple{c}\tuple{a}^{-1} \tuple{b} \in \mathrm{N}^{-}$ by \cref{eq:block_mult} as 
		$\tuple{w} \in \mathrm{N}^{\completion}_{I, \mathbf{r}_I}(m,n)$.\par
		We now show that $\sigma_{\overline{\tuple{x}}}$ is surjective and its fibres have the required cardinality. More precisely 
		we will show that, for all $\tuple{z} \in  \mathrm{N}^{-}$,
			\[
				\sigma_{\overline{\tuple{x}}}^{-1}(\tuple{z}) = \left\{\left. \tuple{y} +	\begin{bmatrix}
																	0 	& 	0\\
																	0	&	\tuple{z}
																	\end{bmatrix}
																\,
																\right\rvert
																\,
																\tuple{y} \in \Mat{m \times n}{m - i_\lastindex}{\completion_N} \text{ lifting }\overline{P}\overline{\tuple{x}} \overline{Q}
																 \right\}.
			\]
		We start by proving that the set on the right-hand side is contained in the fibre of $\tuple{z}$. Let $\tuple{z} \in \mathrm{N}^{-}$. For  
		$\tuple{y} \in \Mat{m \times n}{m - i_\lastindex}{\completion_N}$ lifting $\overline{P}\overline{\tuple{x}} \overline{Q}$,
		we set
			\begin{equation*}
				\tuple{w} = \tuple{y} + \begin{bmatrix}
										0 	& 	0\\
										0	&	\tuple{z}
								\end{bmatrix}.
			\end{equation*}
		We show that $\tuple{w} \in P\mathrm{N}_{\overline{\tuple{x}}}Q$ and $\sigma_{\overline{\tuple{x}}}(\tuple{w}) = \tuple{z}$. 
		Indeed, since $\tuple{y}$ lifts $\overline{P}\overline{\tuple{x}} \overline{Q}$, its first principal $(m - i_\ell)$-minor 
		is invertible. Hence, by \cref{eq:block_mult}, there exist 
			\[
				\begin{aligned}
					\tuple{a} 	&\in \GL{m - i_\lastindex}{\completion_N}{},\\					\tuple{b}	&\in \Mat{(m - i_\lastindex) \times (n - m + i_\lastindex)}{}{\completion_N},\\
					\tuple{c}	&\in \Mat{i_\lastindex \times (m - i_\lastindex)}{}{\completion_N},
				\end{aligned}	
			\qquad\qquad\text{such that }
				\begin{aligned}
					\tuple{y} &= 	\begin{bNiceMatrix}
								\tuple{a}	&	\tuple{b}\\
								\tuple{c}	&	\tuple{c}\tuple{a}^{-1}\tuple{b}
							\end{bNiceMatrix}.
				\end{aligned} 		
			\]
		Therefore, it is enough to show that $\tuple{w} \in P\mathrm{N}_{\overline{\tuple{x}}}Q$, as in that case the definition 
		of $\sigma_{\overline{\tuple{x}}}$ will imply that $\sigma_{\overline{\tuple{x}}}(\tuple{w}) = \tuple{z}$. To this end, 
		note that $P\mathrm{N}_{\overline{\tuple{x}}}Q$ is exactly the set of those matrices in 
		$\mathrm{N}^{\completion}_{I, \mathbf{r}_I}(m,n)$ lifting $\overline{P} \overline{\tuple{x}} \overline{Q}$. Clearly, 
		on the one hand, $P\mathrm{N}_{\overline{\tuple{x}}}Q\subseteq \mathrm{N}^{\completion}_{I, \mathbf{r}_I}(m,n)$ 
		and its elements lift $\overline{P} \overline{\tuple{x}} \overline{Q}$. On the other 
		hand, if $\tuple{v} \in \mathrm{N}^{\completion}_{I, \mathbf{r}_I}(m,n)$ and lifts $\overline{P} \overline{\tuple{x}} \overline{Q}$, 
		then $P^{-1} \tuple{v} Q^{-1}$ lifts  $\overline{\tuple{x}}$. Thus $P^{-1} \tuple{v} Q^{-1} \in \mathrm{N}_{\overline{\tuple{x}}}$.\par
		The fact that $\tuple{w}$ lifts $\overline{P} \overline{\tuple{x}} \overline{Q}$ is immediate 
		from its definition. Thus, by the above, we only need to show that $\tuple{w}$ is in $ \mathrm{N}^{\completion}_{I, \mathbf{r}_I}(m,n)$. 
		However this is also clear, because 
			\begin{multline*}
				\begin{bNiceMatrix}
				\mathrm{Id}_{m - i_\lastindex}	&	0\\
				-\tuple{c}\tuple{a}^{-1}					&	\mathrm{Id}_{i_\lastindex}
				\end{bNiceMatrix}
				\tuple{w}
				\begin{bNiceMatrix}
				\mathrm{Id}_{m - i_\lastindex}					& -\tuple{a}^{-1}\tuple{b} \\
				0	& \mathrm{Id}_{n - m + i_\lastindex}
				\end{bNiceMatrix}															=\\
								\begin{bNiceMatrix}
				\mathrm{Id}_{m - i_\lastindex}	&	0\\
				-\tuple{c}\tuple{a}^{-1}					&	\mathrm{Id}_{i_\lastindex}
				\end{bNiceMatrix}
				\left( \tuple{y}  + \begin{bmatrix}
										0 	& 	0\\
										0	&	\tuple{z}
								\end{bmatrix} \right)
				\begin{bNiceMatrix}
				\mathrm{Id}_{m - i_\lastindex}					& -\tuple{a}^{-1}\tuple{b} \\
				0	& \mathrm{Id}_{n - m + i_\lastindex}
				\end{bNiceMatrix}
																						=\\
				\begin{bNiceMatrix}
					\tuple{a}	&	0\\
					0		&	0
				\end{bNiceMatrix}		+	 \begin{bmatrix}
											0 	& 	0\\
											0	&	\tuple{z}
										\end{bmatrix}
				=  \begin{bmatrix}
											\tuple{a} 	& 	0\\
											0	&	\tuple{z}
										\end{bmatrix}.
			\end{multline*}\par
		We finish the proof by showing that 
		the reverse inclusion also holds. Let $\tuple{w} \in \sigma_{\overline{\tuple{x}}}^{-1}(\tuple{z})$, and 
		let 
			\[
				\begin{aligned}
					\tuple{a} 	&\in \GL{m - i_\lastindex}{\completion_N}{},					& \tuple{b}	&\in \Mat{(m - i_\lastindex) \times (n - m + i_\lastindex)}{}{\completion_N},\\
					\tuple{c}	&\in \Mat{i_\lastindex \times (m - i_\lastindex)}{}{\completion_N},	&\tuple{d}	&\in \Mat{i_\lastindex \times (n - m + i_\lastindex)}{}{\completion_N},
				\end{aligned}	
			\;\; \text{ such that }
				\begin{aligned}
					\tuple{w} &= 	\begin{bNiceMatrix}
								\tuple{a}	&	\tuple{b}\\
								\tuple{c}	&	\tuple{d}
							\end{bNiceMatrix}.
				\end{aligned} 		
			\]
		Clearly 
			\[
				\tuple{w} = 	\begin{bNiceMatrix}
								\tuple{a}	&	\tuple{b}\\
								\tuple{c}	&	\tuple{d}
							\end{bNiceMatrix}
							=
							\begin{bNiceMatrix}
								\tuple{a}	&	\tuple{b}\\
								\tuple{c}	&	\tuple{d} - \tuple{z}
							\end{bNiceMatrix}
							+
							\begin{bmatrix}
								0 	& 	0\\
								0	&	\tuple{z}
							\end{bmatrix}.
			\]
		Since, by definition, $\sigma_{\overline{\tuple{x}}}(\tuple{w}) = \tuple{d} - \tuple{c}\tuple{a}^{-1} \tuple{b}$, 
			\[
				\tuple{y} = \begin{bNiceMatrix}
							\tuple{a}	&	\tuple{b}\\
							\tuple{c}	&	\tuple{d} - \tuple{z}
						\end{bNiceMatrix}
				=
				\begin{bNiceMatrix}
					\tuple{a}	&	\tuple{b}\\
					\tuple{c}	&	\tuple{c}\tuple{a}^{-1} \tuple{b}
				\end{bNiceMatrix},
			\]
		which is in $\Mat{m \times n}{m - i_\lastindex}{\completion_N}$ by \cref{eq:block_mult}. 
		Moreover $\tuple{y}$ lifts $\overline{P}\overline{\tuple{x}} \overline{Q}$, because 
		$\tuple{w}$ is a lift of $\overline{P}\overline{\tuple{x}} \overline{Q}$ 
		and $\tuple{z}$ reduces to the zero matrix modulo $\primeideal$. 
		This finishes the proof that $\sigma_{\overline{\tuple{x}}}$ is surjective and that, by Lemma~\ref{lem:lift_rank}, 
		its fibres have cardinality $q^{(N - 1) (m - i_{\ell}) (n + i_{\ell})}$.
	\end{proof}
\begin{rem}
\label{rem:inductive_step}
Note that, combining the last lemma with Lemma~\ref{lem:rank}, we obtain that, for 
$I = \{i_1, \dots, i_\ell \}_{<} \subseteq [m -1]_0$ and $\mathbf{r}_I \in \N^{I}$,
	\begin{multline}
	\label{eq:size_N}
		\lvert \mathrm{N}^{\completion}_{I, \mathbf{r}_I}(m,n) \rvert = \\
			{n\choose n - m + i_\ell}_{q^{-1}}(q^{- i_\ell - 1};q^{-1})_{m - i_\ell}
				 \lvert \mathrm{N}^{\completion}_{I^{-}, \mathbf{r}_I^{-}}(i_\ell, n - m + i_\ell) \rvert
				 q^{N (m - i_{\ell}) (n + i_{\ell})}. 
	\end{multline}
\end{rem}
\begin{prop}
Let $I \subseteq [m - 1]_0$ and $\mathbf{r}_I \in \N^I$. Then 
	\[
		\lvert \mathrm{N}^{\completion}_{I, \mathbf{r}_I}(m,n) \rvert = 
			{n\choose I+n-m}_{q^{-1}} (q^{-(i_1+1)}; q^{-1})_{m-i_1}
			q^{\sum_{i \in I} r_i (m - i ) (n + i)}.
	\]
\end{prop}
	\begin{proof}
		We proceed by induction on the size of $I$. The case $I = \emptyset$ is immediate because, by convention, empty products are 
		equal to $1$ and so the right-hand side of the equality in the statement is $1$. The base of our induction is $\lvert I \rvert = 1$: 
		in this case the equality in hand follows by Lemma~\ref{lem:lift_rank}.\par
		For $\lvert I \rvert > 1$, Remark~\ref{rem:inductive_step} gives the inductive step. 
		Indeed, by the inductive hypothesis,
			\begin{multline*}
				 \lvert \mathrm{N}^{\completion}_{I^{-}, \mathbf{r}_I^{-}}(i_\ell, n - m + i_\ell) \rvert =  \\
				 	{ n - m + i_\ell \choose I^{-} + n - m}_{q^{-1}} ( q^{- i_1 - 1}; q^{-1})_{i_\ell - i_1} 
					q^{\sum_{i \in I^{-}} r_i (i_\ell - i) (n - m + i_\ell + i)}.
			\end{multline*}
		Substituting the last equality in \eqref{eq:size_N} we obtain
			\begin{multline}
			\label{eq:ind_step}
				\lvert \mathrm{N}^{\completion}_{I, \mathbf{r}_I}(m,n) \rvert  = 
					{ n \choose n - m + i_\ell}_{q^{-1}}{ n - m + i_\ell \choose I^{-} + n - m}_{q^{-1}} \cdot\\
					(q^{-i_\ell + 1}; q^{-1})_{m - i_\ell} (q^{- i_1 + 1}; q^{-1})_{i_\ell - i_1} \cdot\\
					q^{N(m - i_\ell) (n + i_\ell) + \sum_{k = 1}^{\ell - 1} r_{i_k} (i_\ell - i_k)(n - m + i_\ell + i_k)}.	
			\end{multline}
		First, the definition of multinomial coefficient implies that
			\[
				{ n \choose n - m + i_\ell}_{q^{-1}}{ n - m + i_\ell \choose I^{-} + n - m}_{q^{-1}} = {n\choose I + n - m}_{q^{-1}}.
			\]
		Secondly, 
			\begin{align*}
				(q^{ - i_\ell - 1} ; q^{-1})_{m - i_\lastindex} 	
									&= \prod_{ k = 0}^{m - i_\ell - 1}(1 - q^{- i_\ell - 1} q^k) \\
									&= \prod_{ k = 0}^{m - i_\ell - 1}(1 - q^{- i_1 - 1} q^{i_1 - i_\ell - k})\\
									&= \prod_{ k = i_\ell - i_1}^{m - i_1 - 1}(1 - q^{- i_1 - 1} q^{- k}).
			\end{align*}
		Therefore
			\[
				(q^{-i_\ell - 1}; q^{-1})_{m - i_\ell} (q^{- i_1 - 1}; q^{-1})_{i_\ell - i_1} 
						= (q^{- i_1 - 1}; q^{-1})_{m - i_1} .
			\]
		Finally, for all $k \in \{1, \dots, \ell - 1\}$, 
			\[
				(m - i_\ell) (n + i_\ell) + (i_\ell - i_k) (n - m + i_\ell + i_k) = (m - i_k) (n + i_k).
			\]
		Thus
			\begin{multline*}
				N(m - i_\ell) (n + i_\ell) + \sum_{k = 1}^{\ell - 1} r_{i_k} (i_\ell - i_k)(n - m + i_\ell + i_k) = \\
				r_{i_\ell} (m - i_\ell) (n + i_\ell) + \sum_{k = 1}^{\ell - 1} r_{i_k} (m - i_k) (n + i_k) = \\
				\sum_{i \in I} r_i (m -i) (n + i).
			\end{multline*}
		The equality in the statement now follows by substituting in \eqref{eq:ind_step}.
	\end{proof}
We finish the computation of $\lvert \mathrm{N}^{\completion}_{I, \mathbf{r}_I}(m,n) \rvert$ by re-writing the multinomial coefficients in a 
more natural way.
\begin{lem}
\label{lemma:01}
Let $a\in\N_0$ and $I = \{ i_1, \dots, i_\ell\}_{<} \subseteq[m - 1]_0$. Then 
\begin{equation*}
{m+a\choose I+a}_X={m\choose I}\frac{(X^{i_1+1+a};X)_{m-i_1}}{(X^{i_1+1};X)_{m-i_1}}.
\end{equation*}
\end{lem}
	\begin{proof}
		For $r,s\in \N_0$ with $r\geq s$, 
		\begin{equation*}
		{r+a\choose s+a}_X={r\choose s}_X\prod_{i=1}^{r-s}\frac{1-X^{s+a+i}}{1-X^{s+i}}.
		\end{equation*}
		Indeed,
		\begin{align*}
		{r+a\choose s+a}_X&=\prod_{i=1}^{r-s}\frac{1-X^{s+a+i}}{1-X^{i}}\\
		&=\prod_{i=1}^{r}\frac{1-X^{i}}{1-X^{i}}\prod_{i=1}^{r-s}\frac{1-X^{s+a+i}}{1-X^{i}}\\
		&= \frac{\prod_{i=1}^{r}(1-X^{i})}{\prod_{i=1}^{r - s}(1-X^{i}) \prod_{i=1}^{s}(1-X^{i})}
				\frac{\prod_{i=1}^{r - s} (1 - X^{s+a+i})}{\prod_{i= s + 1}^{r} (1 - X^{i})}\\
		&={r\choose s}_X\prod_{i=1}^{r-s}\frac{1-X^{s+a+i}}{1-X^{s+i} }.
		\end{align*}
		Thus 
		\begin{align*}
		{m + a\choose I+a}_X&={m + a\choose i_\ell+a}_X{i_\ell+a\choose i_{\ell - 1}+a}_X\cdots{i_2+a\choose i_1+a}_X\\
		&={m\choose I}_X\prod_{j=1}^{m -i_\ell}\frac{1-X^{i_\ell+a+j}}{1-X^{i_\ell+j}}\cdots\prod_{j=1}^{i_2-i_1}\frac{1-X^{i_1+a+j}}{1-X^{i_1+j}}\\
		&={m \choose I}_X\prod_{j=i_1+1}^m\frac{1-X^{j+a}}{1-X^{j}}.
		\end{align*}
	\end{proof}
A direct application of the lemma above gives that, for all $I = \{ i_1, \dots, i_\ell\}_{<}  \subseteq [m - 1]_0$,
	\begin{equation}
	\label{eq:alt_form}
		{n\choose n-m+I}_{X}(X^{i_1 + 1};X)_{m-i_1} = {m\choose I}_{X}(X^{n-m}X^{i_1 + 1};X)_{m-i_1}.
	\end{equation}
To finish the computation, recall that we defined
	\begin{align*}
		\delta_{n,m} 		&= n - m,	\\
		\overline{a}_{m,n}(i)	&= (m - i) ( n + i) + 2i + \delta_{n,m}.
	\end{align*}
By substituting in \eqref{eq:size_N} and using the definition of $f_{m, n}^I(X)$, we have that, for all $I \subseteq [m - 1]_0$ and $\mathbf{r}_I \in \N^{I}$,
	\[
		\lvert \mathrm{N}^{\completion}_{I, \mathbf{r}_I}(m,n) \rvert = f_{m, n}^I(q^{-1}) q^{\sum_{ i \in I} r_i (\overline{a}_{m,n}(i) - 2i - \delta_{n,m})}.
	\]
Using \eqref{eq:p_series} we readily deduce that 
 	\begin{multline*}
		\zirrcomp = 
			\\\frac{1}{1 - q^{\overline{a}_{m,n}(m) - s_2}} \sum_{I \subseteq [m - 1]_0}
									\sum_{\mathbf{r}_{I}\in \N^I} f_{m,n}^I(q^{-1}) 
									q^{\sum_{i \in I} r_i (\overline{a}_{m,n}(i)  - ( m - i) s_1 - s_2)}.								
	\end{multline*}
Computing the sum of the geometric series above we obtain the bivariate representation zeta function in \cref{main:bivariate}.
\begin{rem}
\label{rem:BtoA}
By \cite[Proposition~4.11]{lin2019arith}, 
\[
	\zeta_{G_{m, n}(\completion)}(s) = \big((1 - q^{m + n - s_2}) \zirrcomp \big)\big\lvert_{\begin{subarray}{l} s_1 \to s - 2 \\ s_2 \to m + n \end{subarray}}.
\]
Making the substitution in the last equality, we readily obtain that 
\[
		\zeta_{G_{m, n}(\completion)}(s)=\sum_{I\subseteq[m-1]_0}f_{m, n}^I(q^{-1})\prod_{i\in I}\frac{q^{(m-i)(n+i)-
		(m-i)s}}{1-q^{(m-i)(n+i)-(m-i)s}}.
\]
This proves \cref{main:additive_G}. Moreover, by  \cite[Proposition~1.5]{stavol2011nilpotent}, 
\[
	\sum_{I\subseteq[m-1]_0}{m\choose I}_{X^{-1}} (YX^{-i_1-1};X^{-1})_{m-i_1}
						\prod_{i\in I}\gp{(X^iZ)^{m-i}}=\frac{(X^{-m}YZ;X)_m}{(Z;X)_m}.
\]
Applying the last equality to \cref{main:additive_G}, we deduce that
\begin{equation*}
\zeta_{G_{m,n}(\completion)}=\frac{(q^{-s};q)_m}{(q^{n-s};q)_m}.
\end{equation*}
Substituting this last form of $\zeta_{G_{m,n}(\completion)}$ in the Euler product of $\zeta_{G_{m,n}(\ringint)}$, we easily deduce 
\cref{main:global}. 
The corollary on the topological zeta function follows by taking the limit $q \to 1$ of the local representation zeta function (cf.\ \cite[Section~3]{ros2015topological}).
\end{rem}
\newcommand{\amn}{A(\tuple{X})}
\newcommand{\submat}{M}
\section{The conjugacy class bivariate zeta function}
\label{sec:cc}
We now turn to computing the bivariate conjugacy class zeta function in \cref{main:bivariate}. 
Let $\tuple{X} = (X_1, \dots, X_{m + n})$. We shall give a full description of 
the minors of $\amn$. This will allow us to determine the integral in \cref{prop:padic}.
\subsection{The matrix $\amn$}
We start by computing the matrix $\amn$.  Since the assumption $m \leq n$ does not play a role in this computation, we 
temporarily drop it for this and the next subsection. Moreover, adhering to the notation in \cref{def:cmatrix},
we relabel the central generators of $\Lambda_{m,n}$ so that their indices range in $[mn]$ rather 
than in $[m]\times [n]$. Concretely, we write $z_{(i - 1) n + j} = z_{ij}$ for all $i \in [m],\, j\in [n]$. 
With this notation, the non-zero commutators defining $\Lambda_{m, n}$ are
	\[
		[x_i, x_{m + j}] = z_{(i - 1) n + j}\qquad i \in [m],\, j\in [n].
	\]
We perform a case distinction based on the row in the $A$-commutator matrix.
\subsubsection{First $m$ rows}
Let $i \in [m]$. Then, for $ k \in [m n]$, 
	\begin{align*}
		\lambda_{ij}^k				&=	0							&	j &\in [m]\\
		\lambda_{i, m + j}^k 	&= 	\begin{cases}
								1		& k = (i - 1) n + j\\
								0		& \text{ otherwise}
							\end{cases}					&	j &\in [n].
	\end{align*}
Therefore for $k > i n $ or $k \leq (i -1 ) n$ we have $\amn_{ik} = 0$. For all $k$ such that $(i -1) n < k \leq i n$ 
we have only one $j \in [n]$ such that $\lambda_{i, m+j}^k  \neq 0$. Namely, $j = k - (i -1)n$, for which $\lambda_{i, m+j}^k = 1$. 
It follows that the $i$-th row of $\amn$ is
	\[
		\begin{NiceMatrix}
			0	& \cdots	& 0	& X_{m+1}				& X_{m + 2}	& \cdots	& X_{m+n}	& 0	& \cdots	& 0\\
				& 		& 	& \vert					& 			& 		& 			& 	& 		& 	\\
				& 		& \Block{1-2}{\text{column }(i -1 ) n  + 1.} 	& 			& 		& 			& 	& 		& 	
		\end{NiceMatrix}
	\] 
\subsubsection{Last $n$ rows}
Let $j \in [n]$. Then, for $k \in [mn]$,  
	\begin{align*}
		\lambda_{m + j, m + i}^k					&= 0							&i &\in [n]\\
		\lambda_{m + j, i}^k	= -\lambda_{i, m + j}^k	&= 	\begin{cases}
													-1	& k = (i -1)n + j\\
													0 	& \text{ otherwise}
												\end{cases}				&i &\in [m].			
	\end{align*}
It follows that
	\[
		\amn_{m + j, k} = 	\begin{cases}
						- X_{\frac{k - j + n}{n}}	& k \in \{ (i -1) n + j \mid i \in [m]\}\\
						0					& \text{ otherwise.}
					\end{cases}
	\]
Hence the last $n$ rows of $\amn$ are
	\[\underbrace{
	\begin{matrix}
		\underbrace{
		\begin{NiceMatrix}
			-X_1		&0		&\Cdots	&0		\\
			0		&\Ddots	&		&\Vdots	\\
			\Vdots	&		&\Ddots	&0		\\
			0		&\Cdots	&0		&-X_1
		\end{NiceMatrix}}_{n \text{ columns}}
		&
		\begin{NiceMatrix}
			-X_2		&0		&\Cdots	&0		\\
			0		&\Ddots	&		&\Vdots	\\
			\Vdots	&		&\Ddots	&0		\\
			0		&\Cdots	&0		&-X_2
		\end{NiceMatrix}
		&
		\dots
		&
		\begin{NiceMatrix}
			-X_m		&0		&\Cdots	&0		\\
			0		&\Ddots	&		&\Vdots	\\
			\Vdots	&		&\Ddots	&0		\\
			0		&\Cdots	&0		&-X_m
		\end{NiceMatrix}
	\end{matrix}
	}_{m \text{ blocks}}.
	\]
\begin{exmp}
Let $m = 3$ and $n = 2$. Then $\amn$ is the $5 \times 6$ matrix
	\[
		\begin{pmatrix}
			X_4	& X_5	&		&		&		&		\\
				&		& X_4	& X_5	&		&		\\
				&		&		&		& X_4	& X_5	\\
			-X_1	&		&-X_2	&		& -X_3	& 		\\
				&-X_1	&		&-X_2	& 		& -X_3	\\
		\end{pmatrix},
	\]
where the non-specified entries are zero.
\end{exmp}
\subsection{The minors of $\amn$}
In this section we compute the minors of $\amn$. We start by introducing some notation. 
\begin{notation}
For $h \in [m]$, we write $\tuple{X}_h$ for the $(m + n - 1)$-tuple of variables 
obtained from $\tuple{X}$ by removing $X_h$. Namely,
	\[
		\tuple{X}_h = (X_1, \dots, X_{h - 1}, X_{h + 1}, \dots, X_m, X_{m + 1}, \dots, X_{m + n}).
	\]
For $h \in [n]$, we write $\tuple{X}^h$ for the $(m + n - 1)$-tuple of variables obtained from $\tuple{X}$ by removing 
$X_{m + h}$. Namely,
	\[
		\tuple{X}^h = (X_1, \dots, X_m, X_{m + 1}, \dots, X_{m + h - 1}, X_{m + h + 1}, \dots, X_{m + n}).
	\]
In what follows, $A(i,j)$ is the $A$-commutator matrix of $\Lambda_{i, j}$ for $i,j \in \N$. As already defined,
 we write $A = A(m,n)$.
\end{notation}
%
\subsubsection{The rank over the function field}
We start the computation of the minors of $\amn$ by bounding the rank of $\amn$. Recall that $\numfield_\primeideal$ is 
the fraction field of $\completion$.
\begin{prop}
	The rank of $\amn$ over $\numfield_\primeideal(\tuple{X})$ is at most $m + n -1$.
\end{prop}
\begin{proof}
We shall use induction on $m + n$, but we first need to isolate some initial cases. Namely, if $m = 1$, then $m + n - 1 = n$, 
which is the number of columns of $\amn$, so the statement is trivially true. The same is true if $n = 1$ by symmetry.\par
If $m \neq 1 \neq n$ we show that all $(m + n)$-minors of $\amn$ are zero.
We proceed by induction on $m + n$.  The base 
case is $m = n = 2$, for which 
	\[
		\amn = 	\begin{pmatrix}
					X_3	& X_4	&		&		\\
						&		&X_3	& X_4	\\
					-X_1&		& -X_2	&		\\
						&-X_1	&		& -X_2	\\
				\end{pmatrix}.
	\]
It is immediate to check that the only $4$-minor of this $4 \times 4$ matrix is $0$. Thus the statement holds for $m = n = 2$.\par
Now, let $N  \in \N$ and assume that the statement holds for all $m,n$ such that $m + n < N$. We show
that it also holds for $n,m$ such that $m + n = N$. To this end, let $m \neq 1$ and $n \neq 1$ be such 
that $m + n = N$ and let $\submat$ be an $(m + n)$-square submatrix of $\amn$. 
We show that $\det(M) = 0$.
\begin{enumerate}[label={\underline{Case \arabic*}}, wide = 0pt, listparindent = \parindent]
\item \label{rk_case_1}
Assume there is an $i \in [m]$ such that $M$ is formed by selecting 
no more than one column having $-X_i$ as an entry. In this case, after possibly rearranging 
the rows and columns of $M$, we may assume without loss of generality that $M$ is a submatrix of 
	\[
	A_i (\tuple{X}) = 	
	\scalebox{0.8}{$
				\begin{bNiceArray}{ccc|cccc}[margin]
					\Block{8-3}{A(m - 1, n)(\tuple{X}_i)} &&	& \Block{4-4}{ O_{(m - 1) \times n}}&		&			& 			\\
						&		&	\hspace*{1cm}			& 					&			&			&			\\
						&		&						& 					&			&			&			\\
						&		&						& 					&			&			& 			\\
						&		& 						& -X_i				& 0			& \Cdots		& 0			\\
						&		& 		 				& 0					& \Ddots		&			& \Vdots		\\
						&		& 		 				& \Vdots				&			& \Ddots		& 0			\\
						&		& 		 				& 0					& \Cdots		& \Cdots		& - X_i		\\
					\hline
					0	& \Cdots	& 0						&X_{m + 1}			&\Cdots		& \Cdots		& X_{m + n}
				\end{bNiceArray}
	$}
	\]
and that it is formed by choosing at most one column from the last $n$ columns of $A_i(\tuple{X})$.\par
Note that this is always the case if $m >n$. Indeed the matrix $\amn$ is formed by $m$ groups of $n$ columns, each of the form
	\begin{equation}
	\scalebox{0.8}{$
	\label{eq:col_groups}
		\begin{NiceMatrix}
			0		& 0			&\Block{12-1}{\cdots}	& 0			\\
			\Vdots	& \Vdots		&\hspace*{1cm}		& \Vdots		\\
			0		& 0			&					& 0			\\
			X_{m + 1}	& X_{m + 2}	&					& X_{m + n}	\\
			0		& 0			&					& 0			\\
			\Vdots	& \Vdots		&					& \Vdots		\\
			0		& 			&					& 			\\
			-X_h		& 0			&					& 			\\
			0		& -X_h		&					& 			\\
			\Vdots	& 0			&					& 			\\	
					& \Vdots		&					& 0			\\
			0		& 0			&					& -X_h		
		\end{NiceMatrix},
	$}
	\end{equation}
for $h \in [m]$. Hence the submatrix $M$ cannot be formed by choosing $2$ columns from each of these groups, because 
$2m > m + n$. \par 
The matrix $A(m - 1, n)(\tuple{X}_i)$ has $m + n - 1$ rows. This implies that we are forced to choose the last row of 
 $A_i (\tuple{X})$ when we form an $(m + n)$-square submatrix. Moreover, by assumption, $M$ is formed by choosing 
 at most one of the last $n$ columns of $A_i (\tuple{X})$. Thus, either $\det(M) = 0$ because the last row of $M$ is zero or there 
 are a $j \in [n]$ and some $(n + m - 1)$-square submatrix $M'$ of $A(m - 1, n)(\tuple{X}_i)$ such that
	\[
		\det(M) = \pm X_{m + j} \det(M'),
	\] 
so $\det(M) = 0$ by the inductive hypothesis.\par
\item \label{rk_case_2}
Assume there is a $j \in [n]$ such that $M$ is formed by choosing at most one column, having $X_{m + j}$ as an entry. In this 
case we may assume without loss of generality that $M$ is a submatrix of
	\[
	A^j (\tuple{X}) = 	
	\scalebox{0.8}{$
				\begin{bNiceArray}{ccc|cccc}[margin]
					\Block{8-3}{A(m, n - 1)(\tuple{X}^j)} &&	& X_{m + j}			& 0			& \Cdots		& 0			\\
						&		& 	\hspace*{1cm}			& 0					& \Ddots		&			& \Vdots		\\
						&		& 		 				& \Vdots				&			& \Ddots		& 0			\\
						&		& 		 				& 0					& \Cdots		& \Cdots		& X_{m + j}	\\
						&		&						& \Block{4-4}{ O_{(m - 1) \times n}}&		&			& 			\\
						&		&						& 					&			&			&			\\
						&		&						& 					&			&			&			\\
						&		&						& 					&			&			& 			\\
					\hline
					0	& \Cdots	& 0						&-X_{1}				&\Cdots		& \Cdots		& -X_{m}
				\end{bNiceArray}
	$}
	\]
and $M$ is formed by selecting at most one of the last $m$ columns of $A^j (\tuple{X})$.\par
Note that, as in Case 1, this is the only possible situation when $n > m$ because we may arrange the columns of $\amn$ in $n$ 
groups of the form 
	\begin{equation}
	\label{eq:col_groups_2}
	\scalebox{0.8}{$
		\begin{NiceMatrix}
			X_{m + h}	& 0 			&\Block{12-1}{\cdots}	& 0			\\
			0		& X_{m + h}	&\hspace*{1cm}		& \Vdots		\\
			\Vdots	& 0			&					& 0			\\		
					&\Vdots		&					& X_{m + h}	\\
					& 			&					& 0			\\
					& 			&					& \Vdots		\\	
			\Vdots	& \Vdots		&					&			\\
			0		& 0			&					& 0			\\
			-X_{1}	&-X_2		&					& -X_m		\\
			0		& 0			&					& 0			\\
			\Vdots	& \Vdots		&					& \Vdots		\\
			0		& 0			&					& 0	
		\end{NiceMatrix}
	$}
	,
	\end{equation}
for $h \in [n]$. Since $2n > m + n$, it follows that, when we form $M$, we may not choose two columns for each of these groups.
\par
As in \ref{rk_case_1}, the assumptions imply that the last row of $M$ has at most one non-zero entry. 
Thus either $\det(M) = 0$ because its last row is zero or there are an $i \in [m]$ and some $(n + m - 1)$-square 
submatrix $M'$ of $A(m , n - 1)(\tuple{X}^j)$ such that
	\[
		\det(M) = \pm X_{i} \det(M'),
	\] 
so $\det(M) = 0$ by the inductive hypothesis.\par
\item \label{rk_case_3}
The only remaining case is when $m = n$ and each row of $M$ has more than one non-zero entry. Indeed, assume that
a row of $M$ has only one non-zero entry and originates as row $i$ of $\amn$. Then $M$ must have been formed 
choosing only one column having $- X_i$ as an entry and we are in \ref{rk_case_1}. Similarly,  if a row of $M$ has only one 
non-zero entry and comes from row $m + j$ of $\amn$ for some $j \in [n]$, then we are in \ref{rk_case_2}, because 
$M$ has been formed choosing only one column of $\amn$ having $X_{m+j}$ as an entry.\par
The last paragraph implies that, forming $M$, we have to choose at least two columns from the $m$ groups in \eqref{eq:col_groups}. In fact, since 
we are defining a $2m$-square submatrix, we may choose at most two columns for each of those groups. As a consequence, there are exactly two 
non-zero entries for each of the first $m$ rows of $M$. An analogous argument using the groups of columns in \eqref{eq:col_groups_2} yields that 
there are also exactly two non-zero entries for each of the last $m$ rows. Thus, there are 
$M_{11}, M_{12}, M_{21}, M_{22} \in \Mat{m \times m}{}{\Z[\tuple{X}]}$, 
such that $M$ is equivalent to	
	\begin{multline*}	
		\begin{bNiceArray}{c|c}[margin]
			M_{11}	& M_{12}	\\
			\hline
			M_{21}	& M_{22}	
		\end{bNiceArray} = \\
		\scalebox{0.8}{$
		\begin{bNiceArray}{ccccc|ccccc}[margin]
			X_{m + j_1}	& X_{m + j_2}	& *		& \Cdots	& *		& 0			& 0			& *		& \Cdots	& *		\\
			0			& 0			&\Vdots	&		& \Vdots	& X_{m + j_1}	& X_{m + j_2}	& \Vdots	&		& \Vdots	\\
			\Vdots		& \Vdots		&		&		& 		& 0			& 0			& 		&		&		\\
						& 			&		&		& 		& \Vdots		& \Vdots		& 		&		&		\\
			0			& 0			& *		& \Cdots	& * 		& 0			& 0			& *		& \Cdots	& *		\\
			\hline
			-X_{i_1}		& 0			&\Cdots	&		& 0		& -X_{i_1'}	& 0			&\Cdots	&		& 0		\\
			0			& -X_{i_2}		&		&		& \Vdots	& 0			& -X_{i_2'}	& 		&		& \Vdots	\\
			\Vdots		& 			& \Ddots	&		& 		& \Vdots		& 			& \Ddots	&		&		\\
						& 			&		& \Ddots	& 0		& 			& 			& 		& \Ddots	& 0		\\
			0			& \Cdots		&		& 0		& -X_{i_m}& 0			& \Cdots	 	& 		& 0		& -X_{i_m'}		
		\end{bNiceArray}
		$}
		,
	\end{multline*}
for some $i_1, \dots, i_m, i_1', \dots, i_m' \in [m]$, $j_1, j_2 \in [n]$.
Now, $M_{21}$ and $M_{22}$ commute and $\det(M_{11}) = \det(M_{12}) = 0$. Therefore
	\[
		\det(M) = \pm ( \det(M_{11}) \det(M_{22}) - \det(M_{12}) \det(M_{21})) = 0.
	\]
\end{enumerate}
\end{proof}
\subsubsection{The minors of $\amn$ }
We shall now describe the minors of $\amn$. We introduce the following notation.
\begin{defn}
	Let $m,n \in \N$. For all $k < m + n$ we denote the ideal generated by the $k$-minors of $\amn$ by $\mathcal{J}_{m,n}^k$.
\end{defn}
\begin{defn}
	Let $m,n \in \N$. For all $k < m + n$, we define the set $\mathcal{M}_{m,n}^k(\tuple{X})$ as the set of monomials
		\[
			X_{i_1}X_{i_2}\cdots X_{i_\lambda}\,X_{m + j_1}X_{m + j_2}\cdots X_{m + j_\omega},
		\]
	for some $\lambda\in [n]_0$, $\omega \in [m]_0$ such that $\lambda + \omega = k$, $i_1,\dots, i_\lambda \in [m]$ and 
	$j_1,\dots, j_\omega  \in [n]$.
\end{defn}
In what follows we shall see that, for all $k$, $\mathcal{M}_{m,n}^k(\tuple{X})$ is a generating set for  $\mathcal{J}_{m,n}^k$. We 
start with the following description of the degree of the polynomials in  $\mathcal{J}_{m,n}^k$.
\begin{lem}
	\label{lem:homogeneous}
	Let $m,n \in \N$ and $k \in \N$ such that $k < m + n$. Let $f \neq 0$ be a $k$-minor of $\amn$ and let $M$ be a submatrix of $\amn$ 
	such that $f = \det(M)$. Assume that $M$ is formed by choosing $\omega$ rows from the first $m$ rows of $\amn$ and 
	$\lambda$ rows from the last $n$ rows of $\amn$. Then $f$ is homogeneous of degree $\lambda$ when viewed as a polynomial 
	in $X_1, \dots, X_m$ and is homogeneous of degree $\omega$ when viewed as a polynomial in $X_{m + 1}, \dots, X_{m + n}$.
\end{lem}
	\begin{proof}
		By definition of $\amn$ the last $\lambda$ rows of $M$ have entries in 
		$\Z[X_1, \dots, X_m]$, while the first $\omega$ rows of $M$ have entries in $\Z[X_{m + 1}, \dots, X_{m + n}]$. 
		Since every non-zero
		entry of $M$ is a linear homogeneous polynomial and $\det(M) \neq 0$, the Leibniz formula for the determinant 
		of $M$ shows that all terms of $f$ have degree $\omega$ in $X_1, \dots, X_m$ and degree 
		$\lambda$ in $X_{m + 1}, \dots, X_{m + n}$.
	\end{proof}
\begin{prop}
	Let $m,n \in \N$. Then, for all $k < m + n$, the ideal $\mathcal{J}_{m,n}^k$ is generated by $\mathcal{M}_{m,n}^k(\tuple{X})$.
\end{prop}
\begin{proof}
	We shall prove the following stronger fact. For all $k < m + n$,
		\begin{enumerate}[label=(\emph{\alph*})]
			\item \label{prop:min_1} all non-zero $k$-minors of $\amn$ are in $\mathcal{M}_{m,n}^k(\tuple{X})$ up to sign,
			\item \label{prop:min_2} every monomial in $\mathcal{M}_{m,n}^k(\tuple{X})$ is a minor of $\amn$ up to sign.
		\end{enumerate}
	We proceed by induction on $m + n$. The base case is $m = n = 1$. Clearly, the non-zero minors of 
		\[
			\begin{pNiceMatrix}
				X_2\\
				-X_1
			\end{pNiceMatrix}
		\]
	are $-X_1$ and $X_2$, so both statements \ref{prop:min_1} and \ref{prop:min_2} hold.  We split the inductive step 
	in two subcases.\par
	\begin{enumerate}[label={\underline{Case \arabic*}}, wide = 0pt, listparindent = \parindent]
	\item	\label{case:1}
	Assume that $m \neq 1$ and $m \geq n$ so that $2m \geq m + n > k$. 
	We start by proving \ref{prop:min_1}. Note that, by Lemma~\ref{lem:homogeneous}, it suffices to show that every 
	non-zero minor of $\amn$ is a monomial with coefficient $1$ or $-1$.\par
	Let $M$ be a $k$-square submatrix of $\amn$. Since $2m > k$, there is an
	$i \in[m]$ such that $-X_i$ is an entry of at most one column of $M$. Hence 
	we may assume, without loss of generality, that $M$ is a submatrix of 
		\[
			A_i(\tuple{X}) = 
				\scalebox{0.8}{$
				\begin{bNiceArray}{ccc|cccc}[margin]
					\Block{8-3}{A(m - 1, n)(\tuple{X}_i)} &&	& \Block{4-4}{ O_{(m - 1) \times n}}&		&			& 			\\
						&		&	\hspace*{1cm}			& 					&			&			&			\\
						&		&						& 					&			&			&			\\
						&		&						& 					&			&			& 			\\
						&		& 						& -X_i				& 0			& \Cdots		& 0			\\
						&		& 		 				& 0					& \Ddots		&			& \Vdots		\\
						&		& 		 				& \Vdots				&			& \Ddots		& 0			\\
						&		& 		 				& 0					& \Cdots		& \Cdots		& - X_i		\\
					\hline
					0	& \Cdots	& 0						&X_{m + 1}			&\Cdots		& \Cdots		& X_{m + n}
				\end{bNiceArray}
				$}
				,
		\]
	and that it is formed by selecting at most one from the last $m$ columns of $A_i(\tuple{X})$. There are now 
	$4$ subcases:
		\begin{enumerate}
			\item the matrix $M$ is a submatrix of $A(m - 1, n)(\tuple{X}_i)$. In this case, by the inductive hypothesis, 
			its determinant is a monomial with coefficient $\pm1$.
			\item The last row or the last column of $M$ is zero. This gives $\det(M) = 0$.
			\item The last row of $M$ is 
				\[
					\begin{matrix}
						0	& \cdots & 0	& X_{m + j}.
					\end{matrix}
				\]
				Then $\det(M) = \pm X_{m + j} \cdot \det(M') $ where $M'$ is a submatrix of 
				$A(m - 1, n)(\tuple{X}_i)$. Thus, by the inductive hypothesis, $\det(M)$ is a monomial with coefficient 
				$\pm 1$.
			\item The last column of $M$ has only one non-zero entry, namely $-X_i$. Similar to the previous subcase, 
			a Laplace expansion along the last column and the inductive hypothesis give that $\det(M)$ is a monomial 
			with coefficient $\pm 1$.		
			\end{enumerate}
	We now prove \ref{prop:min_2}. Let $f \in \mathcal{M}_{m,n}^k(\tuple{X})$. Then
		\[
			f = X_{i_1}X_{i_2}\cdots X_{i_\lambda}\,X_{m + j_1}X_{m + j_2}\cdots X_{m + j_\omega},
		\]
	for some $\lambda \in [n]_0$, $\omega \in [m]_0$ such that $\lambda + \omega = k$, $i_1,\dots, i_\lambda \in [m]$. We 
	prove that $f = \pm \det(M)$ for some submatrix $M$ of $A_i(\tuple{X})$. Since $2m > k$, there is an 
	$i \in [m]$ such that $X_i^2 \nmid f$. We have three subcases.
		\begin{enumerate}
			\item Assume that $X_i \nmid f$ and there is no $j\in [n]$ such that $X_{m + j} \mid f$. Then 
				\[
					f = X_{i_1}X_{i_2}\cdots X_{i_\lambda}
				\] 
				and $k = \lambda \leq n \leq m + n - 1$. Thus, by the inductive hypothesis, $f$ is -- up to sign -- a minor of 
				$A(m - 1, n)(\tuple{X}_i)$. Since the latter is a submatrix of $A_i(\tuple{X})$, we conclude that 
				$f$ is a minor of $A_i(\tuple{X})$ up to sign.
				\par
			\item Assume that $X_i \nmid f$ and there is $j \in [n]$ such that $X_{m + j} \mid f$. Then 
				$f / X_{m + j}$ is in  $\mathcal{M}_{m - 1 , n}^{k - 1}(\tuple{X}_i)$ up to sign.
				Hence, by the inductive hypothesis, 
				there is a submatrix $M'$ of $A(m - 1, n)(\tuple{X}_i)$ such that 
					\[
						\frac{f}{X_{m + j}} = \pm \det(M').
					\]
				Since $A(m - 1, n)(\tuple{X}_i)$ is a submatrix of $A_i(\tuple{X})$, 
					\[
						M =	\begin{bNiceMatrix}
								\Block{3-3}{M'}	&		&	&*			\\
											&		&	&\Vdots		\\
											&		&	&*			\\
									0		&\Cdots	& 0	& X_{m + j}	\\
							\end{bNiceMatrix}
					\]
				is also submatrix of $A_i(\tuple{X})$. Clearly $f = \pm \det(M)$; therefore, it is a minor of $A_i(\tuple{X})$ up to sign.
			\item Assume that $X_i \mid f$. In this case $f/ X_{i}$ is in  $\mathcal{M}_{m - 1,n}^{k - 1}(\tuple{X}_i)$ up to 
				sign.	Thus, by the inductive hypothesis, there is a submatrix $M'$ of $A(m - 1, n)(\tuple{X}_i)$ such that
					\[
						\frac{f}{X_{i}} = \pm \det(M').
					\]
				Moreover, the degree of $f/X_{i}$ in $X_1, \dots, X_m$ is at most $n - 1$. Thus by Lemma~\ref{lem:homogeneous}, 
				not all of the last $n$ rows of $A(m - 1, n)(\tuple{X}_i)$ appear in $M'$. 
				Hence $A_i(\tuple{X})$ has a submatrix 
					\[
						M =	\begin{bNiceMatrix}
								\Block{3-3}{M'}	&		&	&0			\\
											&		&	&\Vdots		\\
											&		&	&0			\\
									*		&\Cdots	& *	& X_{i}		\\
							\end{bNiceMatrix},
					\]
				and $f = \pm \det(M)$.
		\end{enumerate}
	\item Assume that $n \neq 1$ and $m < n$. 
	Then there are two permutation matrices $P$ and $Q$ such that 
		\[
			P \amn Q = A(n,m)(-X_{m + 1}, \dots, - X_{m + n}, -X_1, \dots, -X_m).
		\]
	Thus the case $n >m$ follows from \ref{case:1}, after relabelling the variables.
\end{enumerate}
\end{proof}
\subsection{The bivariate conjugacy class zeta function}
We reinstate the assumption $m \leq n$. The integral from \cref{prop:padic} has integrand 
	\begin{equation}
	\label{eq:integrand}
		\prod_{k =1}^{m + n -1} \frac{\|\mathcal{J}_{m,n}^k (\tuple{X})\cup w\mathcal{J}_{m,n}^{k-1}(\tuple{X}) \|_\primeideal}
				{\| \mathcal{J}_{m,n}^{k - 1} (\tuple{X})\|_\primeideal}.
	\end{equation}
Let $\tuple{x} \in W_{m + n}(\completion)$ and $w \in \primeideal$. We start by computing the value of each factor  $\tuple{x}$. 
We perform a case distinction based on $k$.
\begin{enumerate}[label={\underline{Case \arabic*.}}, wide = 0pt, listparindent = \parindent]
\item  Assume $k \leq \min(m,n) = m$. In this case for all $i \in [m]$ and $j\in [n]$
	\begin{align*}
		X_i^{k}, X_{m + j}^{k} 		&\in \mathcal{J}_{m,n}^k (\tuple{X}),\\
		X_i^{k -  1}, X_{m + j}^{k - 1}	& \in \mathcal{J}_{m,n}^{k - 1} (\tuple{X}).
	\end{align*}
Thus $\|\mathcal{J}_{m,n}^k (\tuple{x})\cup w\mathcal{J}_{m,n}^{k-1}(\tuple{x}) \|_\primeideal = 1$ and 
				$\| \mathcal{J}_{m,n}^{k - 1} (\tuple{x})\|_\primeideal = 1$.
\item  Assume $m < k \leq n$. Let $M = \nu(x_1, \dots, x_m)$, $N = \nu(x_{m + 1}, \dots, x_{m + n})$.
	\begin{enumerate}
		\item $M \leq N$ (thus $M = 0$). In this case for all $i \in [m]$,  
				\[
					X_i^{k}	\in \mathcal{J}_{m,n}^k (\tuple{X}),	\quad 	
							X_i^{k -  1}	 \in \mathcal{J}_{m,n}^{k - 1} (\tuple{X}).
				\]
			Thus $\|\mathcal{J}_{m,n}^k (\tuple{x})\cup w\mathcal{J}_{m,n}^{k-1}(\tuple{x}) \|_\primeideal = 1$ and 
			$\| \mathcal{J}_{m,n}^{k - 1} (\tuple{x})\|_\primeideal = 1$.
		\item $0 = N < M$. Let $j \in [n]$ such that $x_{m + j}$ is invertible and let $i \in [m]$ such that 
			$\nu(x_i) = M$. Then $X_i^{k - m} X_{m + j}^m$ has minimal valuation among the monomials of 
			$\mathcal{J}_{m,n}^k (\tuple{X})$ evaluated at $\tuple{x}$. Hence 
				\[
					\| \mathcal{J}_{m,n}^k (\tuple{X}) \|_\primeideal = q^{-(k - m) M}.
				\]
			Similarly, $\| \mathcal{J}_{m,n}^{k-1} (\tuple{X}) \|_\primeideal = q^{-(k - m - 1) M}$. It follows that, 
			for all $w \in \primeideal$, 
				\begin{align*}
					\| \mathcal{J}_{m,n}^{k - 1} (\tuple{x})\|_\primeideal 	&= \max(q^{-(k - m) M}, \| w\|_\primeideal q^{-(k - m - 1) M})\\
															&= q^{-(k - m - 1) M} \max(q^{-M}, \| w \|)\\
															&= q^{-(k - m - 1) M} \| x_1, \dots, x_m, w\|_\primeideal.
				\end{align*}
			Hence
				\[
					 \frac{\|\mathcal{J}_{m,n}^k (\tuple{X})\cup w\mathcal{J}_{m,n}^{k-1}(\tuple{X}) \|_\primeideal}
				{\| \mathcal{J}_{m,n}^{k - 1} (\tuple{X})\|_\primeideal} = \| x_1, \dots, x_m, w\|_\primeideal.
				\]
	\end{enumerate}
\item Assume $k > n$. We set $M = \nu(x_1, \dots, x_m)$ and $N = \nu(x_{m + 1}, \dots, x_{m + n})$.
	\begin{enumerate}
		\item $M \leq N$. Let $i \in [m]$ such that $x_i$ is invertible and let $j \in [n]$ such that $\nu(x_{m + j}) = N$. Since $k >n$, 
		the monomial $X_i^n X_{m + j}^{k - n}$ has minimal valuation among the monomials in 
		$\mathcal{J}_{m,n}^k (\tuple{X})$. Thus 
			\[
				\| \mathcal{J}_{m,n}^{k - 1} (\tuple{x})\|_\primeideal = 
						q^{-(k - n - 1) N} \| x_{m + 1}, \dots, x_{m + n}, w\|_\primeideal
			\] 
		and $\| \mathcal{J}_{m,n}^{k-1} (\tuple{X}) \|_\primeideal = q^{-(k - n - 1) N}$. Hence 
				\[
					 \frac{\|\mathcal{J}_{m,n}^k (\tuple{X})\cup w\mathcal{J}_{m,n}^{k-1}(\tuple{X}) \|_\primeideal}
				{\| \mathcal{J}_{m,n}^{k - 1} (\tuple{X})\|_\primeideal} = \| x_{m + 1}, \dots, x_{m + n}, w\|_\primeideal.
				\]
		\item $N < M$. A similar argument as the one for $M \leq N$ shows that 	
				\[
					 \frac{\|\mathcal{J}_{m,n}^k (\tuple{X})\cup w\mathcal{J}_{m,n}^{k-1}(\tuple{X}) \|_\primeideal}
				{\| \mathcal{J}_{m,n}^{k - 1} (\tuple{X})\|_\primeideal} = \| x_{1}, \dots, x_{m}, w\|_\primeideal.
				\]
	\end{enumerate}
We summarise the computations above in the following table. 

	\begin{table}[h]
		\caption{Values of the $k$-th factor in the product \eqref{eq:integrand} at $\tuple{x}$.}
		\label{tab:integrand}
		\begin{tabular}{lcc}
		\toprule
		$k$ 
		respect to 
		$m$ and $n$		&	\multicolumn{2}{c}{$k$-th factor}							  						\\
		\cmidrule(lr){2-3}
						&	Case $0 = M \leq N$						& Case $0 = N < M$							\\ 
		\midrule
		$1 \leq k\leq m$	&		1								&		1								\\
		$m < k \leq n$	&		1								&	$\| x_1, \dots, x_m, w\|_\primeideal$			\\
		$n < k < m + n$		& $\| x_{m + 1}, \dots, x_{m + n}, w\|_\primeideal$	&	$\| x_1, \dots, x_m, w\|_\primeideal$			\\
		\bottomrule
		\end{tabular}
	\end{table}
\par
From this we deduce that 
	\begin{equation}
		\prod_{k =1}^{m + n -1} \frac{\|\mathcal{J}_{m,n}^k (\tuple{x})\cup w\mathcal{J}_{m,n}^{k-1}(\tuple{x}) \|_\primeideal}
				{\| \mathcal{J}_{m,n}^{k - 1} (\tuple{x})\|_\primeideal} 
					=	\begin{cases}
							\| x_{m + 1}, \dots, x_{m + n}, w\|_\primeideal^{m - 1}	&	M \leq N	\\
							\| x_{1}, \dots, x_{m}, w\|_\primeideal^{n - 1}		&	N < M.	\\
						\end{cases}
	\end{equation}
\end{enumerate}
\newcommand{\za}{Z_1}
\newcommand{\zb}{Z_2}
\newcommand{\zc}{Z_3}
\subsubsection{Computation of the integral}
Splitting the domain of integration according to Table~\ref{tab:integrand}, 
	\begin{equation}
	\label{eq:split_domain}
		\zcc = (1 - q^{mn - s_2})^{-1} ( 1 +  \za + \zb + \zc),
	\end{equation}
where
	\begin{align*}
		\za	\hspace{-3pt}& = \hspace{-4pt}														
		\int_{(w,\tuple{x}) \in \primeideal \times W_m(\completion)  \times W_n(\completion)} 
															\lvert w \rvert _\primeideal^{(m + n - 1) s_1 + s_2 - mn -  2} d\mu\\
		\zb	\hspace{-3pt}& = \hspace{-4pt}														
			 \int_{(w,\tuple{x}) \in \primeideal \times W_m(\completion) \times \primeideal^{(n)}} 
			  \hspace{-3pt}																	
			\lvert w \rvert _\primeideal^{(m + n - 1) s_1 + s_2 - mn -  2} \| x_{m + 1}, \dots, x_{m + n}, w\|_\primeideal^{ - (m - 1) (s_1 + 1)}d\mu\\
		\zc	\hspace{-3pt}& = \hspace{-4pt}														
		\int_{(w,\tuple{x}) \in \primeideal \times \primeideal^{(m)}  \times W_n(\completion)} \lvert w \rvert _\primeideal^{(m + n - 1) s_1 + s_2 - mn -  2} 
			\| x_{1}, \dots, x_{m}, w\|_\primeideal^{ - (n - 1) (s_1 + 1)} d\mu.
	\end{align*}
By \cite[Proposition~2.2]{lin2019bivFGH}, 
	\begin{align*}
		\za = &\big(1 - q^{-m}\big)\big(1 - q^{-n}\big)\frac{\big(1 - q^{-1}\big) q^{-(m + n -1)s_1 - s_2 + mn + 1}}{\big(1 - q^{-(m + n -1)s_1 - s_2 + mn + 1}\big)}\\
		\zb = &\big(1 - q^{-m}\big) \frac{\big(1 - q^{-1}\big) \big(1 - q^{-(m + n - 1)s_1 - s_2 + (m - 1)n + 1}\big) q^{-ns_1 - s_2 + (m - 1) (n + 1) + 1}}
			{\big(1 - q^{-(m + n -1)s_1 - s_2 + mn + 1}\big) \big(1 - q^{-ns_1 - s_2 + (m - 1) (n + 1) + 1}\big)}\\
		\zc = &\big(1 - q^{-n}\big) \frac{\big(1 - q^{-1}\big) \big(1 - q^{-(m + n - 1)s_1 - s_2 + (n - 1)m + 1}\big) q^{-ms_1 - s_2 + (n - 1) (m + 1) + 1}}
			{\big(1 - q^{-(m + n -1)s_1 - s_2 + mn + 1}\big) \big(1 - q^{-ms_1 - s_2 + (n - 1) (m + 1) + 1}\big)}.
	\end{align*}
Substituting in \cref{eq:split_domain} finishes the proof of the second part of \Cref{main:bivariate}. \Cref{main:class_number} follows from
equation (1.2) in \cite{lin2019arith}, that is 
	\[
		\zeta_{G_{m,n}(\completion)}^{\ccn}(s)  = \mathcal{Z}_{G_{m,n}(\completion)}^{\mathrm{cc}}(0, s).
	\]
This last equality, together with \cite[(1.3)]{lin2019arith} also gives that 
	\[
		\mathcal{Z}_{G_{m,n}(\completion)}^{\mathrm{cc}}(0, s) = \mathcal{Z}_{G_{m,n}(\completion)}^{\mathrm{irr}}(0, s), 
	\]
which holds for arbitrary $q$ and infinitely many values of $s$. Thus, equating the formula in \cref{main:class_number} 
and that in \cref{cor:weyl}, gives the formal identity in \cref{cor:joint}. 
\bibliography{../../A/Biblio/Database}

\begin{bibdiv}
\begin{biblist}

\bib{avn2011rationalabs}{article}{
      author={Avni, Nir},
       title={Arithmetic groups have rational representation growth},
        date={2011},
        ISSN={0003-486X},
     journal={Ann. of Math. (2)},
      volume={174},
      number={2},
       pages={1009\ndash 1056},
         url={http://dx.doi.org/10.4007/annals.2011.174.2.6},
      review={\MR{2831112}},
}

\bib{berderonnpir2013uniform}{article}{
      author={Berman, Mark~N.},
      author={Derakhshan, Jamshid},
      author={Onn, Uri},
      author={Paajanen, Pirita},
       title={Uniform cell decomposition with applications to {C}hevalley
  groups},
        date={2013},
        ISSN={0024-6107},
     journal={J. Lond. Math. Soc. (2)},
      volume={87},
      number={2},
       pages={586\ndash 606},
         url={http://dx.doi.org/10.1112/jlms/jds056},
      review={\MR{3046287}},
}

\bib{carschvol2017traceless}{article}{
      author={Carnevale, Angela},
      author={Shechter, Shai},
      author={Voll, Christopher},
       title={Enumerating traceless matrices over compact discrete valuation
  rings},
        date={2018},
        ISSN={0021-2172},
     journal={Israel J. Math.},
      volume={227},
      number={2},
       pages={957\ndash 986},
         url={https://doi.org/10.1007/s11856-018-1755-4},
      review={\MR{3846349}},
}

\bib{dunvol2017uniform}{article}{
      author={Dung, Duong~H.},
      author={Voll, Christopher},
       title={Uniform analytic properties of representation zeta functions of
  finitely generated nilpotent groups},
        date={2017},
        ISSN={0002-9947},
     journal={Trans. Amer. Math. Soc.},
      volume={369},
      number={9},
       pages={6327\ndash 6349},
         url={https://doi.org/10.1090/tran/6879},
      review={\MR{3660223}},
}

\bib{hrumar2018definable}{article}{
      author={Hrushovski, Ehud},
      author={Martin, Ben},
      author={Rideau, Silvain},
       title={Definable equivalence relations and zeta functions of groups
  (with an appendix by {R}af {C}luckers)},
        date={2018},
        ISSN={1435-9855},
     journal={J. Eur. Math. Soc. (JEMS)},
      volume={20},
      number={10},
       pages={2467\ndash 2537},
         url={https://doi.org/10.4171/JEMS/817},
      review={\MR{3852185}},
}

\bib{zap2006zeta}{article}{
      author={{Jaikin-Zapirain}, Andrei},
       title={Zeta function of representations of compact {$p$}-adic analytic
  groups},
        date={2006},
        ISSN={0894-0347},
     journal={J. Amer. Math. Soc.},
      volume={19},
      number={1},
       pages={91\ndash 118 (electronic)},
         url={http://dx.doi.org/10.1090/S0894-0347-05-00501-1},
      review={\MR{2169043 (2006f:20029)}},
}

\bib{laktho1994counting}{article}{
      author={Laksov, Dan},
      author={Thorup, Anders},
       title={Counting matrices with coordinates in finite fields and of fixed
  rank},
        date={1994},
        ISSN={0025-5521},
     journal={Math. Scand.},
      volume={74},
      number={1},
       pages={19\ndash 33},
         url={https://doi.org/10.7146/math.scand.a-12477},
      review={\MR{1277786}},
}

\bib{lubmag1985varieties}{article}{
      author={Lubotzky, Alexander},
      author={Magid, Andy~R.},
       title={Varieties of representations of finitely generated groups},
        date={1985},
        ISSN={0065-9266},
     journal={Mem. Amer. Math. Soc.},
      volume={58},
      number={336},
       pages={xi+117},
         url={https://doi.org/10.1090/memo/0336},
      review={\MR{818915}},
}

\bib{lin2019arith}{article}{
      author={Macedo Lins~de Araujo, Paula},
       title={Bivariate representation and conjugacy class zeta functions
  associated to unipotent group schemes, {I}: {A}rithmetic properties},
        date={2019},
        ISSN={1433-5883},
     journal={J. Group Theory},
      volume={22},
      number={4},
       pages={741\ndash 774},
         url={https://doi.org/10.1515/jgth-2018-0115},
      review={\MR{3975690}},
}

\bib{lin2019bivFGH}{article}{
      author={Macedo Lins~de Araujo, Paula},
       title={Bivariate representation and conjugacy class zeta functions
  associated to unipotent group schemes, {II}: {G}roups of type {$F$}, {$G$},
  and {$H$}},
        date={2020},
        ISSN={0218-1967},
     journal={Internat. J. Algebra Comput.},
      volume={30},
      number={5},
       pages={931\ndash 975},
         url={https://doi.org/10.1142/S0218196720500265},
      review={\MR{4135016}},
}

\bib{brivol2015enumerating}{article}{
      author={O'Brien, E.~A.},
      author={Voll, C.},
       title={Enumerating classes and characters of {$p$}-groups},
        date={2015},
        ISSN={0002-9947},
     journal={Trans. Amer. Math. Soc.},
      volume={367},
      number={11},
       pages={7775\ndash 7796},
         url={https://doi.org/10.1090/tran/6276},
      review={\MR{3391899}},
}

\bib{plarap1994alggroup}{book}{
      author={Platonov, Vladimir},
      author={Rapinchuk, Andrei},
       title={Algebraic groups and number theory},
      series={Pure and Applied Mathematics},
   publisher={Academic Press, Inc., Boston, MA},
        date={1994},
      volume={139},
        ISBN={0-12-558180-7},
        note={Translated from the 1991 Russian original by Rachel Rowen},
      review={\MR{1278263}},
}

\bib{ros2015topological}{article}{
      author={Rossmann, Tobias},
       title={Topological representation zeta functions of unipotent groups},
        date={2016},
        ISSN={0021-8693},
     journal={J. Algebra},
      volume={448},
       pages={210\ndash 237},
         url={http://dx.doi.org/10.1016/j.jalgebra.2015.09.050},
      review={\MR{3438311}},
}

\bib{ros2018average_i}{article}{
      author={Rossmann, Tobias},
       title={The average size of the kernel of a matrix and orbits of linear
  groups},
        date={2018},
        ISSN={0024-6115},
     journal={Proc. Lond. Math. Soc. (3)},
      volume={117},
      number={3},
       pages={574\ndash 616},
         url={https://doi.org/10.1112/plms.12159},
      review={\MR{3857694}},
}

\bib{ros2020average_ii}{article}{
      author={Rossmann, Tobias},
       title={The average size of the kernel of a matrix and orbits of linear
  groups, {II}: duality},
        date={2020},
        ISSN={0022-4049},
     journal={J. Pure Appl. Algebra},
      volume={224},
      number={4},
       pages={106203, 28},
         url={https://doi.org/10.1016/j.jpaa.2019.106203},
      review={\MR{4021917}},
}

\bib{rosvol2019graphs}{article}{
      author={Rossmann, Tobias},
      author={Voll, Christopher},
       title={Groups, graphs, and hypergraphs: average sizes of kernels of
  generic matrices with support constraints},
        date={2019-08},
     journal={Mem. Amer. Math. Soc. (to appear)
  \href{https://arxiv.org/abs/1908.09589}{arXiv:1908.09589}},
}

\bib{sno2012thesis}{thesis}{
      author={Snocken, Robert},
       title={Zeta functions of groups and rings},
        type={Ph.D. Thesis},
     address={\url{https://eprints.soton.ac.uk/372833/}},
        date={2014},
}

\bib{stavol2013hyper}{article}{
      author={Stasinski, Alexander},
      author={Voll, Christopher},
       title={A new statistic on the hyperoctahedral groups},
        date={2013},
     journal={Electron. J. Combin.},
      volume={20},
      number={3},
       pages={Paper 50, 23},
         url={https://doi.org/10.37236/3242},
      review={\MR{3118958}},
}

\bib{stavol2011nilpotent}{article}{
      author={Stasinski, Alexander},
      author={Voll, Christopher},
       title={Representation zeta functions of nilpotent groups and generating
  functions for {W}eyl groups of type {$B$}},
        date={2014},
        ISSN={0002-9327},
     journal={Amer. J. Math.},
      volume={136},
      number={2},
       pages={501\ndash 550},
         url={http://dx.doi.org/10.1353/ajm.2014.0010},
      review={\MR{3188068}},
}

\bib{stazor2018rational}{article}{
      author={Stasinski, Alexander},
      author={Zordan, Michele},
       title={Rationality of representation zeta functions of compact $p$-adic
  analytic groups},
        date={2020-07},
     journal={\href{https://arxiv.org/abs/2007.10694}{arXiv:2007.10694}},
}

\bib{vol2019delta}{article}{
      author={Voll, Christopher},
       title={Ideal zeta functions associated to a family of class-2-nilpotent
  {L}ie rings},
        date={2020},
        ISSN={0033-5606},
     journal={Q. J. Math.},
      volume={71},
      number={3},
       pages={959\ndash 980},
         url={https://doi.org/10.1093/qmathj/haaa010},
      review={\MR{4142717}},
}

\bib{zor2016adjoint}{article}{
      author={Zordan, Michele},
       title={Adjoint orbits of matrix groups over finite quotients of compact
  discrete valuation rings and representation zeta functions},
        date={2018},
        ISSN={0022-2518},
     journal={Indiana Univ. Math. J.},
      volume={67},
       pages={1683\ndash 1709},
}

\end{biblist}
\end{bibdiv}

\end{document}